\newif\ifconfver
\confverfalse     

\newif\ifplainver  
\plainvertrue

\ifplainver
    \confverfalse   
\fi

\ifconfver
     \documentclass[10pt,twocolumn,twoside]{IEEEtran}
\else
    \ifplainver
        \documentclass[11pt]{article}
        \usepackage{fullpage}
    \else
        \documentclass[11pt,draftcls,onecolumn]{IEEEtran}
    \fi
\fi
\usepackage{calc,amsfonts,amssymb,amsmath,bm,url,color,theorem,graphicx,cite,shortcuts_FW}
\usepackage{psfrag,subfigure,float,hyperref}
\usepackage{algorithm}
\usepackage{algorithmic}

\definecolor{orange}{RGB}{255,107,0}

\newcommand{\tcb}[1]{\textcolor{black}{#1}}

\newtheorem{Fact}{Fact}
\newtheorem{Lemma}{Lemma}

\newtheorem{Theorem}{Theorem}
\newtheorem{Assumption}{Assumption}

\newtheorem{Corollary}{Corollary}

\begin{document}

\bibliographystyle{IEEEtran}

\newcommand{\papertitle}
{
Decentralized Frank-Wolfe Algorithm for Convex and Non-convex Problems
}

\newcommand{\paperabstract}{
...
}


\ifplainver


    \title{\papertitle}

    \author{
    Hoi-To Wai, Jean Lafond,  Anna Scaglione, Eric Moulines
    \thanks{The work of H.-T. Wai is supported by NSF CCF-1553746, NSF CCF-1531050 and J. Lafond by
Direction G\'en\'erale de l'Armement and
the labex LMH (ANR-11-LABX-0056-LMH).
Preliminary versions of this work are presented at ICASSP 2016 \cite{defw16} and
GlobalSIP 2016 \cite{globalsip16}.}
        \thanks{H.-T. Wai and A. Scaglione are with the School of Electrical, Computer and Energy Engineering, Arizona State University, Tempe, AZ 85281, USA. E-mails: \texttt{\{htwai,Anna.Scaglione\}@asu.edu}. J. Lafond was with Institut Mines-Telecom, Telecom ParisTech, CNRS LTCI, Paris, France. Email: \texttt{lafond.jean@gmail.com}. E.~Moulines is with
        CMAP, Ecole Polytechnique, Palaiseau, France. Email: \texttt{eric.moulines@polytechnique.edu}.}
    }

    \maketitle

\else
    \title{\papertitle}

    \ifconfver \else {\linespread{1.1} \rm \fi

    \author{Hoi-To Wai, \emph{Student Member, IEEE}, Jean Lafond, Anna Scaglione, \emph{Fellow, IEEE}, and Eric Moulines
    \thanks{The work of H.-T. Wai is supported by NSF CCF-1553746, NSF CCF-1531050 and J. Lafond by
Direction G\'en\'erale de l'Armement and
the labex LMH (ANR-11-LABX-0056-LMH).
Preliminary versions of this work are presented at ICASSP 2016 \cite{defw16} and
GlobalSIP 2016 \cite{globalsip16}.}
        \thanks{H.-T. Wai and A. Scaglione are with the School of Electrical, Computer and Energy Engineering, Arizona State University, Tempe, AZ 85281, USA. E-mails: \texttt{\{htwai,Anna.Scaglione\}@asu.edu}. J. Lafond is with Institut Mines-Telecom, T\'{e}l\'{e}com ParisTech, Universit\'{e} Paris Saclay, Paris, France. Email: \texttt{lafond.jean@gmail.com}. E.~Moulines is with
        CMAP, \'{E}cole Polytechnique, Palaiseau, France. Email: \texttt{eric.moulines@polytechnique.edu}.}
         }

    \maketitle

    \ifconfver \else
        \begin{center} \vspace*{-2\baselineskip}
        \end{center}
    \fi

    \ifconfver \else \IEEEpeerreviewmaketitle} \fi

 \fi
 
     \begin{abstract}
Decentralized optimization algorithms have received much attention due to the
recent advances in
network information processing.
However, conventional decentralized algorithms based on projected gradient descent are
incapable of handling high dimensional constrained problems,
as the projection step becomes computationally prohibitive to compute.
To address this problem, this paper adopts a projection-free optimization approach,
a.k.a.~the Frank-Wolfe (FW) or conditional gradient algorithm.
We first develop a decentralized FW (DeFW) algorithm from the classical
FW algorithm. The convergence of the proposed algorithm is studied by viewing the
decentralized algorithm as an \emph{inexact} FW algorithm.
Using a diminishing step size rule and letting $t$ be the iteration number,
we show that the DeFW algorithm's convergence rate is ${\cal O}(1/t)$
for convex objectives; is ${\cal O}(1/t^2)$ for strongly convex objectives with the optimal solution
in the interior of the constraint set; and is ${\cal O}(1/\sqrt{t})$ towards
a stationary point for smooth but non-convex objectives.
We then show that a {  consensus}-based
DeFW algorithm meets the above guarantees
with two communication rounds per iteration.
Furthermore, we demonstrate the advantages of the proposed DeFW algorithm
on low-complexity robust matrix completion
and communication efficient sparse learning.
Numerical results on synthetic and real data
are presented to support our findings.\end{abstract}

\ifplainver
\else
    \begin{IEEEkeywords}\vspace{-0.0cm}
       decentralized optimization, high-dimensional optimization, Frank-Wolfe algorithm, consensus algorithms, communication efficient algorithms, LASSO, matrix completion
    \end{IEEEkeywords}
\fi

\ifconfver \else
    \ifplainver \else
        \newpage
\fi \fi

\section{Introduction}

\ifplainver
Recently, 
\else
\IEEEPARstart{R}ecently, 
\fi
algorithms for tackling high-dimensional optimizations have been sought
for the 
`big-data' challenge \cite{Cevher2014}.
As these data/measurements are dispersed over clouds of networked machines,
it is important to consider \emph{decentralized} algorithms that
can allow the agents/machines to co-operate and
leverage the aggregated computation power \cite{Sayed2013}.

This paper considers decentralized algorithm for
tackling a constrained optimization problem 
with $N$ agents:
\beq \label{eq:opt}
\min_{ \prm \in \RR^d }~F(\prm)~{\rm s.t.}~\prm \in \Cset,~~\text{where}~~F(\prm) \eqdef \frac{1}{N} \sum_{i=1}^N f_i (\prm) \eqs,
\eeq
and $f_i (\prm )$ is a continuously differentiable (possibly non-convex)
function held by the $i$th agent and 
$\Cset$ is a closed and bounded convex set in $\RR^d$.
Typically, the function $f_i(\prm)$ models the loss function 
over the private data
held by agent $i$, and
$\Cset$ corresponds to a regularization constraint imposed on the 
optimization variable $\prm$ such
as sparsity or low rank of an array of values. 
Problem~\eqref{eq:opt} covers a number of applications in {control
theory},
signal processing and machine learning,
including
{system identification \cite{liu10}},
matrix completion \cite{Candes_Recht09} and
sparse learning \cite{ravazzi15,patterson14}.

To tackle \eqref{eq:opt}, 
the agents communicate on a network described as a graph $G= (V,E)$,
where $V = \{1,...,N\}$ is the set of agents and $E \subseteq V \times V$ describes the connectivity.
As $G$ is not fully connected, it is useful to apply decentralized algorithms
that relies on \emph{near-neighbor} information exchanges.
In light of this, various authors have proposed to tackle \eqref{eq:opt} through
decentralized algorithms that are built on the average
consensus protocol \cite{Tsitsiklis1984,Dimakis2010}. For example,
\cite{bianchi13,ichen12,na16,johan08,RNV2010,Shi2015,jako14,nedich16}
studied the decentralized gradient methods; 
\cite{xli13,varagnolo16} considered the Newton-type methods; \cite{th14,duichi12,wei13,Simonetto2016,mingyi16} 
considered the decentralized
primal-dual algorithms.
\cite{aldo14,aldo16} are built on the successive convex approximation framework.
The convergence properties of these algorithms were investigated extensively,
especially for convex objectives \cite{ichen12,na16,Sayed2013,Simonetto2016,nedich16,johan08,xli13,RNV2010,wei13,mingyi16,Shi2015,th14,jako14,duichi12,varagnolo16};
for non-convex objectives, a few recent results can be found in
\cite{mingyi16,aldo14,aldo16,bianchi13,htwai15,de_dict_learn}.
However, most prior work are \emph{projection}-based such that
each iteration of the above
algorithms 
necessitates a {projection} step 
onto the constraint set $\Cset$, or solving a sub-problem with 
similar complexity.
While the projection step entails a modest complexity
for problems with moderate dimension,
it may become computationally prohibitive for
high dimensional problems,
thereby rendering the 
above methods impractical for our scenario of interest.

To address the issue above, 
this paper focuses on a decentralized \emph{projection-free} algorithm.
We extend the Frank-Wolfe (FW) algorithm \cite{fw56}
to operate with near-neighbor information exchange, leading to a decentralized FW (DeFW) algorithm.
The FW (a.k.a.~conditional gradient) algorithm has been 
recently popularized due to its efficacy in handling
high-dimensional constrained problems.
Examples of its applications include
{optimal control \cite{wu83},}
matrix completion \cite{jaggi10}, image and video colocation \cite{joulin14},
electric vehicle charging \cite{zhang16}
and traffic assignment \cite{transport84}; see the overview article \cite{Jaggi13}.
From the algorithmic perspective, the FW algorithm replaces the costly projection step
in PG based algorithms 
with a constrained linear optimization, which often
admits an efficient solution.
In a centralized setting, the convergence of FW algorithm has been  
studied for convex problems \cite{fw56,Jaggi13}, yet
little is known for non-convex problems \cite{henry15,simon16,ourpaper}. 

Our contributions are as follows. We first describe abstractly the DeFW algorithm
as a variation of the FW algorithm with \emph{inexact} iterates and gradients. 
We then analyze its convergence ---
for convex objectives,
the sub-optimality (of objective values) of the iterates produced by the proposed algorithm
is shown to converge as ${\cal O}(1/t)$ with $t$ being the iteration number, and is
${\cal O}(1 / t^2)$ for strongly convex objectives when the optimal solution
is in the interior of $\Cset$; for non-convex objectives,
we demonstrate that the proposed algorithm has limit points that are stationary points of \eqref{eq:opt},
and they can be found at the rate of ${\cal O}( \sqrt{1/t} )$.
We then show that a {consensus}-based
implementation of the proposed algorithm
with fixed number of communication rounds of near-neighbor information exchanges achieves all the above guarantees.
To our knowledge, this is the first decentralized FW algorithm.
Moreover, our convergence rate in the non-convex setting is comparable
to that of a \emph{centralized} projected gradient method \cite{lan_ncvx_15}.
Lastly, we present examples on communication-efficient LASSO 
and decentralized matrix completion.

The rest of this paper is organized as follows.
Section~\ref{sec:algo} develops the DeFW algorithm.
We summarize the main theoretical results of convergence
for convex and non-convex objective functions.
A consensus-based implementation of the DeFW algorithm
will then be presented in Section~\ref{sec:gossip}.
Applications of the DeFW algorithm are discussed in Section~\ref{sec:ext}.
Finally, Section~\ref{sec:num} presents numerical results to
 support our  findings.

\subsection{Notations \& Mathematical Preliminaries}
For any $d \in \NN$, we define $[d]$ as the set $\{1,...,d\}$.
We use boldfaced lower-case letters to denote vectors and boldfaced upper-case letters to denote matrices.
For a vector ${\bm x}$ (or a matrix ${\bm X}$),
the notation $[ {\bm x} ]_i$ (or $[{\bm X}]_{i,j}$) denotes its $i$th element (or $(i,j)$th element).
The vectorization of a matrix ${\bm X} \in \RR^{m_1 \times m_2}$ is denoted by ${\rm vec} ({\bm X}) = [{\bm x}_1; {\bm x}_2;
\ldots ; {\bm x}_{m_2}] \in \RR^{m_1 m_2}$ such that ${\bm x}_i$ is the $i$th column of ${\bm X}$.
The vector ${\bf e}_i \in \RR^d$ is the $i$th unit vector such that $[ {\bf e}_i ]_j = 0$ for all $j \neq i$ and
$[{\bf e}_i]_i = 1$.
For some positive finite constants $C_1,C_2,C_3$, $C_2 \leq C_3$ 
and non-negative functions $f(t), g(t)$,
the notations $f(t) = {\cal O}( g(t) )$, $f(t) = \Theta ( g(t))$
indicate $f(t) \leq C_1 g(t)$, $C_2 g(t) \leq f(t) \leq C_3 g(t)$ for sufficiently
large $t$, respectively.

Let $\E$ be a Euclidean space
embedded in $\RR^d$ and the Euclidean norm is denoted by $\| \cdot \|_2$.
The binary operator $\langle \cdot, \cdot \rangle$ denotes the inner product on $\E$.
In addition, $\E$ is equipped with a norm $\| \cdot \|$ and the corresponding dual norm $\| \cdot \|_\star$.
Let $G, L ,\mu$ be non-negative constants.
Consider a function $f : \RR^d \rightarrow \RR$, the function $f$ is  $G$-Lipschitz if for
all $\prm, \prm' \in \E$,
\beq
|f( \prm ) - f( \prm' )| \leq G \| \prm - \prm' \| \eqs;
\eeq
the function $f$ is $L$-smooth if for
all $\prm, \prm' \in \E$,
\beq
f( \prm) - f(\prm') \leq \langle \grd f(\prm'), \prm - \prm' \rangle + \frac{L}{2} \| \prm - \prm' \|_2^2 \eqs,
\eeq
the above is equivalent to $\| \grd f(\prm') - \grd f(\prm) \|_2 \leq L \| \prm' - \prm \|_2$;
the function $f$ is $\mu$-strongly convex if for all $\prm, \prm' \in \E$,
\beq
f( \prm) - f(\prm') \leq \langle \grd f(\prm), \prm - \prm' \rangle - \frac{\mu}{2} \| \prm - \prm' \|_2^2 \eqs;
\eeq
moreover, $f$ is convex if the above is satisfied with $\mu=0$. 

Consider Problem \eqref{eq:opt},
its constraint set $\Cset \subseteq \E$ is convex and bounded with the diameter defined as:
\beq \label{eq:rho}
\rho \eqdef \max_{ \prm, \prm' \in \Cset } \| \prm - \prm' \|,~~\bar{\rho} \eqdef \max_{ \prm, \prm' \in \Cset } \| \prm - \prm' \|_2 \eqs,
\eeq
note that $\rho$ is defined with respect to (w.r.t.) the norm $\| \cdot \|$ while $\bar{\rho}$
is defined w.r.t. the Euclidean norm.
When the objective function $F$ is $\mu$-strongly convex with $\mu >0$,
the optimal solution to \eqref{eq:opt}
is unique and denoted by $\prm^\star$, we also define
\beq \label{eq:int}
\delta \eqdef \min_{ {\bm s} \in \partial \Cset } \| {\bm s} - \prm^\star \|_2 \eqs,
\eeq
where $\partial \Cset$ is the boundary set of $\Cset$.
If $\delta > 0$, the solution $\prm^\star$ is in the interior of $\Cset$.

\section{Decentralized Frank-Wolfe (DeFW)} \label{sec:algo}
We develop the decentralized Frank-Wolfe (DeFW) algorithm from the
classical FW algorithm \cite{fw56}.
Let $t \in \NN$ be the iteration number and the initial point
$\prm_0 \in \Cset$ is feasible. 
{Recall the definition $F(\prm) \eqdef (1/N) \sum_{i=1}^N f_i (\prm)$,}
the \emph{centralized} FW algorithm for problem \eqref{eq:opt} proceeds by:
\begin{subequations} \label{eq:fw}
\begin{align}
\label{eq:lo} \atom_{t-1} & \in \arg \min_{\atom \in \Cset}~ \langle \grd F(\prm_{t-1}) , \atom \rangle \eqs, \\
\label{eq:updfw} \prm_{t} & = \prm_{t-1} + \gamma_{t-1} ( \atom_{t-1} - \prm_{t-1} ) \eqs,
\end{align}
\end{subequations}
where $\gamma_{t-1} \in (0,1]$ is a step size to be determined.
Observe that $\prm_{t}$ is a convex combination of $\prm_{t-1}$ and $\atom_{t-1}$ \tcb{which are both feasible}, therefore
$\prm_{t} \in \Cset$ as $\Cset$ is a convex set.
When the step size is chosen as $\gamma_t = 2/(t+1)$, the FW algorithm
is known to converge at a rate of ${\cal O}(1/t)$ if $F$ is $L$-smooth
and convex \cite{Jaggi13}.
A main feature of the FW algorithm is that 
the linear optimization\footnote{Notice that 
\eqref{eq:lo} is a convex optimization problem with a linear objective.
} (LO) \eqref{eq:lo} can be solved more efficiently 
than computing a projection, leading to a \emph{projection-free} algorithm.
At the end of this section, we will 
illustrate a few examples of $\Cset$ with efficient
LO computations.

\algsetup{indent=1em}
\begin{algorithm}[t]
\caption{Decentralized Frank-Wolfe (DeFW).}\label{alg:defw}
  \begin{algorithmic}[1]
  \STATE \textbf{Input}: Initial point $\prm_1^i$ for $i=1,\dots,N$.
  \FOR {$t=1,2,\dots$}
     \STATE \label{line:con} \emph{Consensus}: approximate the average iterate: \label{fw:con}
   \beq \notag
   \bar{\prm}_t^i \leftarrow \texttt{NetAvg}_t^i ( \{ \prm_t^j \}_{j=1}^N ),~~\forall~i \in [N] \eqs. \vspace{-0.4cm}
   \eeq
   \STATE \label{line:agg} \emph{Aggregating}: approximate the average gradient: \label{fw:agg}
   \beq \notag \textstyle
   \bargrd{t}{i} \leftarrow \texttt{NetAvg}_t^i ( \{ \grd f_j( \bar{\prm}_t^j ) \}_{j=1}^N ),~~\forall~i \in [N] \eqs. \vspace{-0.4cm}
   \eeq
   \STATE \label{line:fw} \emph{Frank-Wolfe Step}: update
   \beq \notag \hspace{-0.2cm} \prm_{t+1}^{i} \leftarrow (1-\gamma_t) \bar{\prm}_t^i + \gamma_t \atom_t^i~~\text{where}~~
   \atom_t^i \in \arg \min_{ \atom \in \Cset } \langle \bargrd{t}{i}, \atom \rangle,
   \eeq
   for all agent $i \in [N]$ and $\gamma_t \in (0,1]$ is a step size. \vspace{.1cm}
\ENDFOR
\STATE \textbf{Return}: $\bar{\prm}_{t+1}^i, \forall~i  \in [N]$.
  \end{algorithmic}
\end{algorithm}

Our next endeavor is to extend the FW algorithm to a decentralized setting via mimicking
\eqref{eq:fw} with only near-neighbor information exchanges.
Doing so requires replacing the centralized gradient/iterate
$\grd F(\prm_t)$, $\prm_t$ in \eqref{eq:fw} with  
local approximations, as similar to the strategy in \cite{johan08,Simonetto2016}.

In the following,  we offer a high-level description of the proposed 
DeFW algorithm and discuss the convergence properties of it.
Details regarding the implementation will be postponed to Section~\ref{sec:gossip}.
Let $\prm_t^i$ denotes an auxillary iterate kept by agent $i$ at iteration $t$. Define 
the average iterate:
\beq
\textstyle \bar{\prm}_t
\eqdef N^{-1} \sum_{i=1}^N \prm_t^i \eqs
\eeq
and the local iterate
$\bar{\prm}_t^i$ as an approximation
of the average iterate above, also kept by agent $i$.
We require 
$\bar{\prm}_t^i$
to track 
$\bar{\prm}_t$ with an increasing accuracy.
Let $\{ \dprm_t \}_{t \geq 1}$ be a non-negative, decreasing sequence 
with $\dprm_t \rightarrow 0$, we assume
\begin{Assumption}[\protect{$\{ \dprm_t \}_{t \geq 1}$}] \label{ass:a1}
\tcb{For all $t \geq 1$, it holds that}
\beq \label{eq:dprm}
\textstyle \max_{i \in [N]} \| \bar{\prm}_t^i - \bar{\prm}_t \|_2 \leq  \dprm_t \eqs.
\eeq
\end{Assumption}
To compute \eqref{eq:lo}, ideally each agent
has to access the \emph{global gradient}, $\grd F( \bar{\prm}_t )$.
However, just the local function $f_i(\cdot)$ is available 
and agent $i$ can 
only compute the local gradient
$\grd f_i( \bar{\prm}_t^i )$. 
Therefore, we also need to track 
the average gradient,
\beq
\textstyle \bargrd{t}{} \eqdef N^{-1} \sum_{j=1}^N \grd f_j( \bar{\prm}_t^j )  \eqs,
\eeq
by the local approximation $\bargrd{t}{i}$.
{Note that $\bargrd{t}{}$ is close to $\grd F( \bar{\prm}_t )$
when each of the function $f_i(\prm)$ is smooth and
$\bar{\prm}_t^i$ is close to $\bar{\prm}_t$.}
Let $\{ \dgrd_t \}_{t \geq 1}$ be a non-negative, decreasing sequence
with $\dgrd_t \rightarrow 0$, we assume:
\begin{Assumption}[\protect{$\{ \dgrd_t \}_{t \geq 1}$}] 
\label{ass:a2}
\tcb{For all $t \geq 1$, it holds that}
\beq \label{eq:dgrd} \textstyle
 \max_{i \in [N]} \| \bargrd{t}{i} - \bargrd{t}{} \|_2 \leq \dgrd_t \eqs.
\eeq
\end{Assumption}
Naturally, from
the local approximation $\bargrd{t}{i}$, 
the $i$th agent can compute the update direction
$\atom_t^i = \arg \min_{ \atom \in \Cset } \langle \bargrd{t}{i}, \atom \rangle$ 
and update $\prm_{t+1}^i$ similarly as in \eqref{eq:updfw}.
To summarize, a sketch of the DeFW algorithm
can be found in Algorithm~\ref{alg:defw}.

Under Assumptions \ref{ass:a1}-\ref{ass:a2},
for each agent $i$, line~\ref{line:fw} in Algorithm~\ref{alg:defw} can be
regarded as performing an 
\emph{inexact} FW update on $\bar{\prm}_t$,
whose convergence can be characterized below.
For convex objective functions, we have:
\begin{Theorem} \label{thm:cvx}
Set the step size as $\gamma_t = 2/(t+1)$. Suppose that each of $f_i$ is convex and $L$-smooth.
\tcb{Let $C_p$ and $C_g$ be two positive constants}.
Under Assumptions~\ref{ass:a1}-\ref{ass:a2} 
[$\dprm_t = C_p / t$, $\dgrd_t = C_g / t$], we have
\beq \label{eq:cvx_re1}
F(\bar{\prm}_t) - F(\prm^\star) \leq  \frac{8\bar{\rho}(C_g + L C_p) + 2 L \bar{\rho}^2}{t+1}\eqs,
\eeq
for all $t \geq 1$,
where $\prm^\star$ is an optimal solution to \eqref{eq:opt}.
Furthermore, if $F$ is $\mu$-strongly convex and the optimal solution $\prm^\star$
lies in the interior of $\Cset$, i.e.,  $\delta > 0$ (cf.~\eqref{eq:int}), we have
\beq \label{eq:cvx_re2}
F(\bar{\prm}_t) - F(\prm^\star) \leq \frac{ (4\bar{\rho}(C_g + L C_p) + L \bar{\rho}^2 )^2 }{2 \delta^2 \mu} \cdot \frac{9}{(t+1)^2} \eqs,
\eeq
for all $t \geq 1$.
\end{Theorem}
The proof can be found in Appendix~\ref{pf:cvx}. 
We remark that $\bar\prm_t$ is always feasible.
For strongly convex objective functions,
the conditions \eqref{eq:cvx_re1}, \eqref{eq:cvx_re2} imply that
the sequence $\{ \bar{\prm}_t \}_{t \geq 1}$ converges to
an optimal solution of
\eqref{eq:opt}.
Furthermore, as the consensus error,
$ \max_{i \in [N]} \| \bar{\prm}_t^i - \bar{\prm}_t \|_2$, 
decay to zero (cf.~Assumption~\ref{ass:a1}),
the local iterates $\{ \bar\prm_t^i \}_{t \geq 1}$ share similar convergence 
guarantee as $\{ \bar{\prm}_t \}_{t \geq 1}$. 

For non-convex objective functions, 
we study the convergence of the FW/duality gap:
\beq 
g_t \eqdef \max_{ \prm \in \Cset } \langle \grd F( \bar{\prm}_t) , \bar{\prm}_t - \prm  \rangle \eqs. 
\eeq
From the definition, when $g_t = 0$, the iterate $\bar\prm_t$
will be a stationary point to \eqref{eq:opt}. Thus we may regard
$g_t$ as a measure of the stationarity of the iterate $\bar{\prm}_t$.
Also, define the set of stationary point to \eqref{eq:opt} as:
\beq \label{eq:stationary} 
\Cset^\star = \big\{ \underline{\prm} \in \Cset : \max_{ \prm \in \Cset } ~\langle \grd F( \underline{\prm}) , \underline{\prm} - \prm  \rangle = 0 \big\} \eqs.
\eeq
We consider the following technical assumption: \vspace{-.2cm}
\begin{Assumption} \label{ass:a3}
{The set $\Cset^\star$ is non-empty.} Moreover,
the function $F(\prm)$ takes a finite number of values over $\Cset^\star$, i.e.,
the set $F( \Cset^\star ) = \{ F( \prm ) : \prm \in \Cset^\star \} $ is finite.
\end{Assumption}
Verifying Assumption~\ref{ass:a3} may be hard in practice.
Meanwhile, it is reasonable to assume
that \eqref{eq:opt} has a finite number
of stationary points since the set $\Cset$ is bounded.
In the latter case, Assumption~\ref{ass:a3} is satisfied. 
We now have:

\begin{Theorem} \label{thm:ncvx}
Set the step size as $\gamma_t = 1/t^{\alpha}$ for some $\alpha \in (0,1]$.
Suppose each of $f_i$ is $L$-smooth and
$G$-Lipschitz (possibly non-convex). \tcb{Let $C_p$, $C_g$ be two positive constants}. 
\tcb{Under Assumption~\ref{ass:a1}-\ref{ass:a2} [$\dprm_t = C_p/t^\alpha$,
$\dgrd_t = C_g/t^\alpha$], it holds that}:
\begin{enumerate}
\item for all $T \geq 6$ that are even, if $\alpha \in [0.5,1)$, 
\beq
\label{eq:mainthm}
\begin{split}
\min_{t \in [T/2+1,T] } ~g_t \leq & ~\frac{1}{T^{1-\alpha}} \cdot \frac{1-\alpha}{(1-(2/3)^{1-\alpha})} \cdot \\
& \hspace{-1.4cm} \Big({G {\rho} + ( L \bar{\rho}^2/2 + 2\bar{\rho} ( C_g + L C_p) ) \log2} \Big) \eqs;
\end{split}
\eeq
if $\alpha \in (0,0.5)$, 
\beq
\label{eq:mainthm2}
\begin{split}
\hspace{-.3cm} \min_{t \in [T/2+1,T] } ~g_t \leq &  ~\frac{1}{T^{\alpha}} \cdot \frac{1-\alpha}{(1-(2/3)^{1-\alpha})}~\cdot \\
\hspace{-.3cm} & \hspace{-2.2cm} \Big( G {\rho} + \frac{( L \bar{\rho}^2/2 + 2\bar{\rho} ( C_g + L C_p) ) (1-(1/2)^{1-2\alpha})}{1-2\alpha} \Big) \eqs.
\end{split}
\eeq
\item
additionally, under Assumption~\ref{ass:a3} and $\alpha \in (0.5,1]$,
the sequence of objective values $\{ F( \bar{\prm}_t )\}_{t\geq 1}$ converges,
$\{ \bar{\prm}_t \}_{t \geq 1}$
has limit points and each limit point is in $\Cset^\star$.
\end{enumerate}
\end{Theorem}
The proof can be found in Appendix~\ref{pf:ncvx}.
Note that
setting $\alpha = 0.5$ gives the quickest convergence rate of ${\cal O}(1/\sqrt{T})$.
It is worth mentioning that our results are novel compared to prior work 
on non-convex FW even in a centralized setting ($N=1, \dprm_t=0, \dgrd_t = 0$). 
For instance, \cite{henry15}
requires that the local minimizer is unique; \cite{simon16} gives the
same convergence rate but uses an adaptive step size.
We remark that the local iterates $\{ \bar\prm_t^i \}_{t \geq 1}$ share similar
convergence property as $\{ \bar{\prm}_t \}_{t \geq 1}$ 
due to Assumption~\ref{ass:a1}.

Lastly,  we  
survey some relevant examples of the constraint set $\Cset$
where the LO in the DeFW algorithm 
(cf.~line~\ref{line:fw} in Algorithm~\ref{alg:defw}) can be computed 
efficiently:
\begin{enumerate}
\item When $\Cset$ is the $\ell_1$ ball, $\Cset = \{ \prm \in \RR^d : \| \prm \|_1 \leq R \}$,
\beq \label{eq:l1}
\atom_t^i = -R \cdot {\bf e}_k,~\text{where}~k \in \arg \max_{ j \in [d] }~ \big| [ \bargrd{t}{i} ]_j \big| \eqs.
\eeq
The solution above amounts to finding the coordinate index of $\bargrd{t}{i}$
with the maximum magnitude.
Importantly, this solution is only $1$-sparse. Consequently, 
the $t$th iterate $\bar{\prm}_t$ will be at most
$t N$-sparse.
The worst-case complexity of computing $\atom_t^i$ is ${\cal O}(d)$; in comparison,
the worst-case complexity for the projection into an
$\ell_1$ ball is ${\cal O}(d \log d )$\footnote{There exists a randomized, accelerated algorithm for projection in \cite{Duchi_Shalev-Shwartz_Singer_Chandra08} with an \emph{expected} complexity of ${\cal O}(d)$.}.
\item When $\Cset$ is the trace norm ball, 
$\Cset = \{ \prm \in \RR^{m_1 \times m_2} : \| \prm \|_{\sigma, 1} \leq R \}$,
where $ \| \prm \|_{\sigma, 1}$ is the sum of the singular values of $\prm$.
Let ${\bm u}_1, {\bm v}_1$ be the top-1 left/right singular vector of $\bargrd{t}{i}$, we have
\beq \label{eq:tracenorm}
\atom_t^i = -R \cdot {\bm u}_1 {\bm v}_1^\top \eqs.
\eeq
Importantly, at a target solution accuracy of $\delta$,
the top singular vectors can be computed with a complexity of
${\cal O} ( \max\{ m_1, m_2 \} \log ( 1 / \delta) )$ 
using the power/Lanczos method
if $\| {\rm vec}( \bargrd{t}{i} ) \|_0 = {\cal O} ( \max\{m_1,m_2\} )$.
In comparison, the projection onto the trace norm ball
requires
a complexity of ${\cal O} ( \max\{m_1 m_2^2 ,m_2 m_1^2\} 
\log ( 1 / \delta)  )$ for computing the full SVD of an
$m_1 \times m_2$ matrix \cite{Golub_VanLoan13}.
\end{enumerate}

For more examples of $\Cset$
with efficient LO computations, the interested readers
are referred to
\cite{Jaggi13} for an overview.

\section{Consensus-based DeFW algorithm} \label{sec:gossip}
This section demonstrates how an average consensus (AC) based scheme
can generate local approximations $\bar{\prm}_t^i$, $\bargrd{t}{i}$ 
at the desirable sub-linear 
accuracies (cf.~Assumptions~\ref{ass:a1}-\ref{ass:a2}) 
using as few as two 
communication rounds per iteration.

Specifically, the following discussions are based on the static AC
\cite{Tsitsiklis1984,Dimakis2010}. 
{While the exact details are left for a future
work, we believe that it is possible to
extend our protocol to a time varying network's setting}.
We consider an undirected graph $G = (V,E)$ and assign
a non-negative, symmetric \emph{weighted adjacency matrix} 
${\bm W} \in \RR_+^{N \times N}$
that describes the local communication between the $N$ agents. The matrix satisfies 
$W_{ij} \eqdef [{\bm W}]_{ij} > 0$ iff $(i,j) \in E$
and it is doubly
stochastic, i.e., ${\bm W} {\bf 1} = {\bm W}^\top {\bf 1} = {\bf 1}$. Moreover,
\begin{Assumption} \label{ass:w}
The second largest  (in magnitude) eigenvalue of ${\bm W}$ is
strictly less than one, i.e., $|\lambda_2( {\bm W})| < 1$.
\end{Assumption} 
The existence of such matrix ${\bm W}$ is guaranteed
if $G$ is connected.
For each round of the AC update, the agents take a weighted average of
the values from its neighbors according to ${\bm W}$.
We now state the following fact regarding ${\bm W}$.
\begin{Fact} \label{fact:gac}
Let ${\bm x}_1, ..., {\bm x}_N \in \RR^d$ be a set of vectors and
${\bm x}_{\sf avg} \eqdef N^{-1} \sum_{i=1}^N {\bm x}_i$ be their average.
Suppose ${\bm W}$ is a doubly stochastic, non-negative matrix.
The output after performing one round of AC update:
\beq \label{eq:gac_prototype} \textstyle
\overline{\bm x}_i = \sum_{j=1}^N W_{ij} \cdot {\bm x}_j
\eeq
must satisfy
\beq \label{eq:gac_presult}
\sqrt{ \sum_{i=1}^N \| \overline{\bm x}_i - {\bm x}_{\sf avg} \|_2^2 } \leq |\lambda_2 ({\bm W})| \cdot \sqrt{\sum_{i=1}^N  \| {\bm x}_i - {\bm x}_{\sf avg} \|_2^2} \eqs,
\eeq
where $\lambda_2 ( {\bm W})$ is the second largest eigenvalue of ${\bm W}$.
\end{Fact}
The fact above can be verified from linear algebra.
Together with Assumption~\ref{ass:w}, the above implies that each AC update
\eqref{eq:gac_prototype} moves
the vectors closer to the average ${\bm x}_{\sf avg}$.
Repeatedly applying \eqref{eq:gac_presult} shows the well known fact that
AC computes the average ${\bm x}_{\sf avg}$ at a linear rate.

Let us consider the near-neighbor computation of $\bar{\prm}_t^i$ in line~\ref{line:con} of 
Algorithm~\ref{alg:defw}. Here, the \emph{consensus step}
is computed by:
\beq \textstyle \label{eq:gac2}
\bar{\prm}_t^i = \sum_{j =1}^N W_{ij} \cdot \prm_t^j \eqs,
\eeq
i.e., we apply one round of the AC update.
Since $W_{ij} = 0$ if $(i,j) \notin E$, the above operation is achievable by
information exchanges with the near-neighbors of agent $i$.

Now, for some $\alpha \in (0,1]$,
we define $t_0(\alpha)$ as the smallest integer such that
\beq \label{eq:t0}
\lambda_2 ({\bm W}) \leq \Big( \frac{ t_0 (\alpha) } {t_0(\alpha) + 1} \Big)^\alpha \cdot \frac{ 1 } { 1 + (t_0(\alpha))^{-\alpha} } \eqs.
\eeq
Notice that $t_0 (\alpha)$ is upper bounded by:
\beq
t_0 (\alpha) \leq \lceil ( |\lambda_2 ({\bm W})|^{-1/(1+\alpha)} - 1 )^{-1} \rceil \eqs,
\eeq
which is finite under Assumption~\ref{ass:w}.
The following lemma can be
easily proven:
\begin{Lemma} \label{lem:gac2}
Set the step size $\gamma_t = 1 / t^\alpha$ in the DeFW algorithm for some $\alpha \in (0,1]$, then $\bar{\prm}_t^i$ in \eqref{eq:gac2} satisfies 
Assumption~\ref{ass:a1}:
\beq \label{eq:gac2r}
\textstyle \max_{i \in [N]} \| \bar{\prm}_t^i - \bar{\prm}_t \|_2 \leq \dprm_t = C_p / t^\alpha,~\forall~t \geq 1 \eqs,
\eeq
\beq
C_p \eqdef (t_0(\alpha))^\alpha \cdot \sqrt{N} \bar{\rho} \eqs.
\eeq
\end{Lemma}
The proof is postponed to Appendix~\ref{pf:gac2}, which relies on the observation
that $\prm_t^i$ is a convex 
combination of $\bar{\prm}_{t-1}^i$ and $\atom_{t-1}^i$. In particular,
$\bar{\prm}_{t-1}^i$ is already ${\cal O}(1/(t-1)^\alpha)$-close to the network average
from the last iteration and the update direction $\atom_{t-1}^i$
is weighted by a decaying step size $\gamma_{t-1}$.

In comparison to what we were able to establish above,
the near-neighbor computation of $\bargrd{t}{i}$ in line~\ref{line:agg} of the DeFW algorithm
is less straightforward.
Unlike the computation of $\bar{\prm}_t$,
computing $N^{-1} \sum_{i=1}^N \grd f_i ( \bar{\prm}_t^i )$ to an accuracy of ${\cal O}(1/t^{\alpha})$
by only communicating
the local gradient $\grd f_i(\bar{\prm}_t^i)$  
requires $\sim \log t $ rounds of
updates when the AC protocol is employed.
The main issue is that the local gradient $\grd f_i ( \bar{\prm}_t^i )$ computed
by the $i$th agent is different from the local gradient computed at the other agent,
even when $\bar{\prm}_t^i$ is close to $\bar{\prm}_t^j$ for $j \neq i$.

We propose an approach that is 
inspired by the fast stochastic average gradient (SAGA)
method \cite{saga14} which re-uses the gradient approximation
$\bargrd{t-1}{i}$ from the last iteration\footnote{After the submission, 
the authors 
notice that a similar technique is adopted
in \cite{na16,nedich16,aldo16} under the name of `gradient tracking'
for various decentralized methods.
We provide a rate analysis with non-asymptotic constants.}.
Define the following surrogate of local gradient at iteration $t$:
\beq \label{eq:grad_svrg}
\bgrd{t}{i} \eqdef \bargrd{t-1}{i} + \grd f_i( \bar{\prm}_t^i ) - \grd f_i(\bar{\prm}_{t-1}^i ) \eqs,
\eeq
for all $ i \in [N]$.
When $t=1$, we set $\bgrd{1}{i} = \grd f_i ( \bar{\prm}_1^i )$.
We now apply the AC update to the gradient surrogate. 
In line~\ref{line:agg} of Algorithm~\ref{alg:defw}, the \emph{aggregating} step
is computed by:
\beq \textstyle \label{eq:gac1}
\bargrd{t}{i} = \sum_{ j = 1 }^N W_{ij} \cdot \bgrd{t}{j} \eqs,
\eeq
i.e., using just one round of AC update on $\bgrd{t}{i}$.
Below we show that the average gradient is preserved by $\bgrd{t}{i}$.
Moreover, $\bargrd{t}{i}$
has an approximation error similar to Lemma~\ref{lem:gac2}:
\begin{Lemma} \label{lem:gac1}
Set the step size $\gamma_t =1 / t^\alpha$ in the DeFW algorithm for some $\alpha \in (0,1]$.
Suppose that each of $f_i$ is $L$-smooth, and $\bar{\prm}_t^i$ is updated according to \eqref{eq:gac2},
then $\bargrd{t}{i}$ in \eqref{eq:gac1} satisfies
\beq \textstyle \label{eq:avg1}
N^{-1} \sum_{i=1}^N \bargrd{t}{i} = N^{-1} \sum_{i=1}^N \bgrd{t}{i} = N^{-1} \sum_{i=1}^N \grd f_i (\bar{\prm}_t^i ),~\forall~t \geq 1 \eqs,
\eeq
and Assumption~\ref{ass:a2}:
\beq \label{eq:gac1r} \textstyle
 \max_{i \in [N]} \| \bargrd{t}{i} - \bargrd{t}{} \|_2 \leq \dgrd_t = C_g / t^\alpha,~\forall~t \geq 1 \eqs,
\eeq
\beq
\begin{split}
& C_g \eqdef \sqrt{N} \max \Big\{ \!~ 2 (2C_p + \bar{\rho} ) L, (t_0(\alpha) )^\alpha |\lambda_2({\bm W})|  \Big( \frac{ L \bar{\rho} }{1 - |\lambda_2({\bm W})|}  + B_1 \Big) \Big\} \eqs,
\end{split}
\eeq
where $B_1 \eqdef \max_{ i=1,...,N } \| \grd f_i ( \bar{\prm}_1^i ) \|_2$ is a bound
on the initial gradients.
\end{Lemma}
The proof can be found in Appendix~\ref{pf:gac1}. 
Similar intuition as in Lemma~\ref{lem:gac2}
was used in the proof. In particular, we observe that $\bargrd{t}{t}$ is
a linear combination of $\bargrd{t-1}{i}$ and 
$\grd f_i( \bar{\prm}_t^i ) - \grd f_i(\bar{\prm}_{t-1}^i )$.
The former is ${\cal O}(1/(t-1)^\alpha)$-close to the average,
and the latter decays to zero as enforced by the step size $\gamma_t$. 

Finally, the conditions on $\dprm_t, \dgrd_t$ necessitated by Theorem~\ref{thm:cvx}
\& \ref{thm:ncvx} can be satisfied by \eqref{eq:gac2} \& \eqref{eq:gac1}.
\begin{Corollary}
Under Assumption~\ref{ass:w},
the results in
Theorem~\ref{thm:cvx} \& \ref{thm:ncvx} hold when line~\ref{line:con}, line~\ref{line:agg}
in the DeFW algorithm (Algorithm~\ref{alg:defw}) 
are computed by  \eqref{eq:gac2}, \eqref{eq:gac1} respectively.
\end{Corollary}
In other words, the consensus-based DeFW algorithm converges for both
convex and non-convex problems, while using a constant number of communication rounds
per iteration.

As a final remark,
recall from Theorem~\ref{thm:ncvx} that for non-convex objectives,
the best rate of convergence
can be achieved if we set $\alpha = 0.5$. However,
from Lemma~\ref{lem:gac2} \& \ref{lem:gac1}, we notice that the approximation error
also decays the slowest when $\alpha = 0.5$. 
This presents a potential tradeoff in the choice of $\alpha$. 
From our numerical experience, we find that
setting $\alpha = 0.75$ yields a
good performance in practice.

\section{Applications} \label{sec:ext}
In this section, we study two applications of the DeFW algorithm to illustrate
its flexibility and efficacy.

\subsection{Example I: Decentralized Matrix Completion} \label{sec:mc}
Consider a setting when the network of agents obtain incomplete observations of a matrix
${\prm}_{\sf true}$ of dimension $m_1 \times m_2$ with $m_1, m_2 \gg 0$.
The $i$th agent has corrupted observations from the \emph{training} set
$\Omega_i \subset [m_1] \times [m_2]$ that
are expressed as:
\beq \label{eq:mat_obs}
Y_{k,l} = [{\prm}_{\sf true}]_{k,l} + Z_{k,l},~\forall~(k,l) \in \Omega_i \eqs.
\eeq
To recover a low-rank $\prm_{\sf true}$, we consider the following
trace-norm constrained matrix completion (MC) problem:
\beq \label{eq:mc}
\min_{ \prm \in \RR^{m_1 \times m_2}}~\sum_{i=1}^N
\sum_{ (k,l) \in \Omega_i } \tilde{f}_i ( [\prm]_{k,l} ,  Y_{k,l} )
~{\rm s.t.}~\| \prm \|_{\sigma,1} \leq R \eqs,
\eeq
where $\tilde{f}_i : \RR^2 \rightarrow \RR$ is a loss function picked by agent $i$
according to the observations he/she received.
{  Notice that \eqref{eq:mc} is also related to the low rank subspace
system identification problem described in \cite{liu10}, where ${\bm Y}$ with
$[{\bm Y}]_{k,l} = Y_{k,l}$,
$\prm_{\sf true}$ are modeled as
the measured system response and the ground truth
low rank response; also see \cite{anna08} for a related work.}

Similar MC problems have been
considered in \cite{qing12,mackey2015distributed,yu2012scalable,recht2013parallel},
where \cite{qing12} studied a consensus-based optimization method similar
to ours and \cite{mackey2015distributed,yu2012scalable,recht2013parallel}
studied the parallel computation setting where the agents are working synchronously
in a fully connected network.
Compared to our approach, these work assume that the rank of $\prm_{\sf true}$
is known in advance and solve the MC problem via matrix factorization.
In addition, \cite{qing12,mackey2015distributed} required that each local observation
set $\Omega_i$ only have entries  taken
from a disjoint subset of the columns/rows only.
Our approach does not have any restrictions above.

We consider two different observation models. When $Z_{k,l}$ is the i.i.d.~Gaussian
noise of variance $\sigma_i^2$,
we choose $\tilde{f}_i ( \cdot , \cdot )$ to be the square loss function,
i.e.,
\beq
\tilde{f}_i ( [\prm]_{k,l},  Y_{k,l} ) \eqdef (1/\sigma_i^2) \cdot ( Y_{k,l} - [ \prm ]_{k,l} )^2 \eqs.
\eeq
This yields the classical MC problem in \cite{Candes_Recht09}.
The next model considers the sparse+low rank matrix completion
in \cite{chandrasekaran2009sparse}, where
the observations are contaminated with a sparse noise. Here, we model
$Z_{k,l}$ as a \emph{sparse} noise in the sense that there are
a few number of entries in $\Omega_i$ where $Z_{k,l}$ is non-zero.
We choose $\tilde{f}_i ( \cdot,\cdot )$ to be the negated Gaussian loss, i.e.,
\beq
 \tilde{f}_i ( [\prm]_{k,l} ,  Y_{k,l} ) \eqdef \Big( 1- {\rm exp} \Big( - \frac{ ( [\prm]_{k,l} - Y_{k,l} )^2 }{\sigma_i} \Big) \Big) \eqs,
\eeq
where $\sigma_i > 0$ controls the robustness to outliers for the data obtained at the $i$th agent.
Here, $\tilde{f}_i(\cdot,\cdot )$ is a \emph{smoothed} $\ell_0$ loss \cite{mohimani07}
with enhanced robustness to outliers in the data.
Notice that the resultant MC problem \eqref{eq:mc}
is non-convex.

Note that \eqref{eq:mc} is a special case of problem \eqref{eq:opt} with
$\Cset$ being the trace-norm ball.
The consensus-based DeFW algorithm can be applied  on \eqref{eq:mc} directly.
The projection-free nature of the DeFW algorithm
leads to a low complexity implementation \eqref{eq:mc}.
Lastly, several remarks on the communication and storage cost 
of the DeFW algorithm are in order:
\begin{itemize}
\item The SAGA-like gradient surrogate $\grd_t^i F$ \eqref{eq:grad_svrg} is supported
only on $\cup_{i=1}^N \Omega_i$ since for all $ i \in [N]$, the local gradient
\beq
\grd f_i( \bar{\prm}_t^i ) = \sum_{(k,l) \in \Omega_i} \tilde{f}_i' ( [ \bar{\prm}_t^i ]_{k,l}, Y_{k,l} ) \cdot {\bf e}_k ({\bf e}_l')^\top
\eeq
is supported on $\Omega_i$, where $\bar{\prm}_t^i$ is defined in \eqref{eq:gac2}.
In the above, ${\bf e}_k$ (${\bf e}_l'$) is the $k$th ($l$th) canonical basis vector for $\RR^{m_1}$ ($\RR^{m_2}$)
and $\tilde{f}_i'(\theta, y)$ is the derivative of $\tilde{f}_i ( \theta, y)$ taken with respect to $\theta$.
Consequently, the average $\bargrd{t}{i}$ is supported only on $\cup_{i=1}^N \Omega_i$.
As $| \cup_{i=1}^N \Omega_i | \ll m_1 m_2$, the amount of information exchanged during the \emph{aggregating}
step (Line~\ref{fw:agg} in DeFW) is low.
\item The update direction $\atom_t^i$ is a rank-one matrix composed of
the top singular vectors of $\bargrd{t}{i}$ (cf.~\eqref{eq:tracenorm}). Since every iteration
in DeFW adds at most $N$ distinct pair of singular vectors to $\bar{\prm}_t$,
the rank of $\bar{\prm}_t^i$ is upper bounded by $t N$ if we
initialize by $\bar{\prm}_0^i = {\bm 0}$. We can reduce the communication
cost in Line~\ref{fw:con} in DeFW by exchanging these singular vectors. Note that $(tN) \cdot (m_1 + m_2)$ entries are stored/exchanged instead of $m_1 \cdot m_2$.

\item
When the agents are \emph{only} concerned with predicting the entries of ${\prm}_{\sf true}$ in the subset
$\Xi \subset [m_1] \times [m_2]$, instead of propagating the singular vectors as
described above, the \emph{consensus}
step can be carried out by exchanging only the entries of $\prm_{t+1}^i$ in
$\Xi \cup \big(\cup_{i=1}^N \Omega_i \big)$ without affecting the operations of
the DeFW algorithm.
In this case, the storage/communication cost is
$|\Xi \cup \big(\cup_{i=1}^N \Omega_i \big)|$.
\end{itemize}

\subsection{Example II: Communication Efficient DeFW for LASSO} \label{sec:lasso}
Let $( {\bm y}_i, {\bm A}_i )$ be the available data tuple at agent $i \in [N]$ such that
${\bm A}_i \in \RR^{m \times d}$ and ${\bm y}_i \in \RR^m$. The data ${\bm y}_i$ is a
corrupted measurement of an unknown parameter $\prm_{\sf true}$:
\beq  \label{eq:sparse}
{\bm y}_i = {\bm A}_i \prm_{\sf true} + {\bm z}_i \eqs,
\eeq
where ${\bm z}_i \sim {\cal N} ( {\bm 0}, \sigma^2 {\bm I} )$ are independent noise vectors.
Furthermore, we assume $m \ll d$ such that the matrix ${\bm A}_i^\top {\bm A}_i$ is rank-deficient.
However, the parameter $\prm_{\sf true}$ is $s$-sparse such that $s = \| \prm_{\sf true} \|_0 \ll d$.
This motivates us to consider the following distributed LASSO problem:
\beq \label{eq:lasso}
\min_{ \prm \in \RR^d } ~\sum_{i=1}^N \frac{1}{2} \| {\bm y}_i -{\bm A}_i \prm \|_2^2 ~~{\rm s.t.}~~\| \prm \|_1 \leq R \eqs,
\eeq
Notice that the above is a special case of \eqref{eq:opt} with
$f_i (\prm) = (1/2) \| {\bm y}_i -{\bm A}_i \prm \|_2^2 $ and $\Cset$ is an 
$\ell_1$-ball in $\RR^d$ with radius $R$.
We assume that \eqref{eq:lasso} has an optimal solution $\prm^\star$ that is sparse.
{ 
The settings above also correspond to identifying a linear system
described by a sparse parameter $\prm_{\sf true}$, where
${\bm A}_i$, ${\bm y}_i$ are the input, output of the system,
respectively; see \cite{bako11} for a related formulation on the
identification of switched linear systems.
}

A number of decentralized algorithms are easily applicable
to \eqref{eq:lasso}.
For example, the decentralized projected gradient 
(DPG) algorithm in \cite{RNV2010} is described by --- at iteration $t$,
\beq\label{eq:dpg} \textstyle
\prm^{i,PG}_{t+1} = {\cal P}_{\Cset} \big( \sum_{j=1}^N W_{ij} \prm^{j,PG}_t - \alpha_t \grd f_i \big( \sum_{j=1}^N W_{ij} \prm^{j,PG}_t \big) \big) \eqs,
\eeq
where $\alpha_t \in (0,1]$ is a diminishing step size and $W_{ij}$ is the weighted
adjacency matrix
described in Section~\ref{sec:gossip}.
For convex problems, the DPG algorithm  is shown to converge
to an optimal solution $\prm^\star$ of \eqref{eq:lasso} at a rate of ${\cal O}(1/\sqrt{t})$ \cite{ichen12}.

Let us focus on  
the communication efficiency of the DPG algorithm, which is important when
the network between agents is limited in bandwidth.
To this end, we define the communication cost
as the \emph{number
of non-zero real numbers exchanged per agent}.
As seen from \eqref{eq:dpg}, at each iteration the $i$th agent exchanges
its current iterate $\prm^{i,PG}_t$ with the neighboring agents.
From the computation step shown, 
$\prm^{i,PG}_t$ may contain as high as ${\cal O}(d)$ non-zeros and
the per-iteration communication cost
will be  ${\cal O} (d)$. 
Despite the high communication cost,
the per-iteration computation complexity of \eqref{eq:dpg}
is also high, i.e., at ${\cal O}(d \log d)$ \cite{Duchi_Shalev-Shwartz_Singer_Chandra08}.
We notice that \cite{ravazzi15,patterson14} have considered distributed sparse recovery
algorithm with focus on the communication efficiency. However, their algorithms are
based on the iterative hard thresholding (IHT) formulation \cite{blumensath08} that requires a-priori knowledge
on the sparsity level of $\prm_{\sf true}$.
Our consensus-based DeFW algorithm in 
Section~\ref{sec:gossip} may also be applied directly
to \eqref{eq:lasso}. 
However, similar issue as the DPG algorithm may arise
during the \emph{aggregating} step, 
since the gradient surrogate \eqref{eq:grad_svrg}
may also have ${\cal O}(d)$ non-zeros.
Lastly, another related work is \cite{yildiz08} which
applies coding to `compress' the message exchanged
in the consensus-based algorithms.

This section proposes a \emph{sparsified DeFW algorithm}
for solving \eqref{eq:lasso}.
The modified algorithm applies a novel 
`sparsification' procedure to 
reduce communication
cost during the iterations,
which is enabled by the structure of the DeFW algorithm.
To describe the sparsified DeFW algorithm, we first argue that
the \emph{consensus step} in the consensus-based DeFW should
remain unchanged
as it already has a low communication cost.
From \eqref{eq:l1} and \eqref{eq:gac2}, we see
that ${\prm}_t^i$ is at most $(t-1)N + 1$-sparse
since $\atom_t^i$ is always a $1$-sparse vector\footnote{As pointed out by \cite{Jaggi13},
this observation also leads to an interesting sparsity-accuracy trade-off
when applying FW on $\ell_1$ constrained problems.} (cf.~\eqref{eq:l1}).
As such, the communication cost of this step is bounded
by $t N$.

Our focus is to improve the communication efficiency of \emph{aggregating step}.
Here, the key idea is that only the \emph{largest magnitude
coordinate} in $\bargrd{t}{i}$ is sought when computing $\atom_t^i$ 
(cf.~\eqref{eq:l1}).
As long as the largest magnitude coordinate in $\bargrd{t}{i}$ is preserved,
the updates in the DeFW algorithm can remain unaffected.
This motivates us to `sparsify' the gradient information
at each iteration before exchanging them with the neighboring agents.
Let $\Omega_t \subseteq [d]$ be the coordinates
of the gradient information to be exchanged at iteration $t$. The agents
exchange the following gradient surrogate in lieu of \eqref{eq:grad_svrg}:
\beq \label{eq:sparse_grad}
\widehat{\grd_t^i F} \eqdef \big(\grd f_i( \bar{\prm}_t^i ) \big) \odot {\bf 1}_{\Omega_t},~~\text{where}~~\textstyle {\bf 1}_{\Omega_t} = \sum_{k \in \Omega_t} {\bf e}_k \eqs,
\eeq
and $\odot$ denotes the Hadamard/element-wise product.

Let $\ell_t = \lceil C_l + \log ( t ) / \log |\lambda_2^{-1} ( {\bm W} ) | \rceil$ where $C_l$
is some finite constant and $\lambda_2 ( {\bm W} )$ is the second largest
eigenvalue of the weight matrix ${\bm W}$,
the sparsified DeFW algorithm computes the approximate gradient average
$\bargrd{t}{i}$ in line~\ref{line:agg} of
Algorithm~\ref{alg:defw} by:
\beq \label{eq:gac3} \textstyle
\bargrd{t}{i} = \sum_{j=1}^N [ {\bm W}^{\ell_t} ]_{ij} \cdot \widehat{\grd_t^j F} \eqs.
\eeq
Note that \eqref{eq:gac3} requires $\ell_t$ \emph{rounds} of AC updates
to be performed at iteration $t$, i.e., a logarithmically increasing number of rounds of
AC updates.
The update direction $\atom_t^i$ can then be computed by sorting 
the vector $\bargrd{t}{i}$.
As $\bargrd{t}{i}$ is $| \Omega_t |$-sparse, this update direction
can be computed in ${\cal O}( |\Omega_t|)$ time.

We pick the coordinate set $\Omega_t$ in a decentralized manner. Consider the following decomposition:
\beq \textstyle
\Omega_t = \bigcup_{i=1}^N \Omega_{t,i} \eqs,
\eeq
where $\Omega_{t,i} \subset [d]$ is picked by agent $i$ at iteration $t$.
The coordinate set $\Omega_t$ needs to be known by all agents before
\eqref{eq:gac3}. This can be achieved with low communication overhead,
e.g., by forming a spanning tree on the graph $G$ and
broadcasting the required indices in $\Omega_t$ to all agents; see \cite{distcomp04}.
Set $p_t$ as the maximum desirable cardinality of $\Omega_{t,i}$,
agent $i$ chooses the coordinate set using one of the following two schemes:
\begin{itemize}
\item \emph{(Random coordinate)} Each agent selects $\Omega_{t,i}$
by picking $p_t$ coordinates uniformly (with replacement) from $[d]$.
\item \emph{(Extreme coordinate)} Each agent selects $\Omega_{t,i}$ as
the $p_t$ largest magnitude coordinates of the vector $\grd f_i (\bar{\prm}_t^i)$.
\end{itemize}

For the random coordinate selection scheme,
the following lemma shows that the gradient approximation error
can be controlled at a desirable rate with an appropriate choice of $p_t$.
Let $\xi_t \eqdef (1 - (1 - 1 /d)^{p_t N} )$, we have:
\begin{Lemma} \label{lem:prob}
Set $\epsilon > 0$ and $\ell_t = \lceil C_l + \log ( t ) / \log |\lambda_2^{-1} ( {\bm W} )| \rceil$.
Let {$p_t \geq C_0 t$ for some $C_0 < \infty$}.
With probability at least $1- \pi^2 \epsilon / 6$, the following holds
for all $\prm \in \Cset$:
\beq
\begin{split}
\hspace{-.1cm} \Big\| \xi_t^{-1} \bargrd{t}{i} - \frac{1}{N} \sum_{i=1}^N \grd f_i (\bar{\prm}_t^i ) \Big\|_\infty &
= {\cal O} \Big( \frac{ d \sqrt{ \log( t^2 / \epsilon ) } }{ t N } \Big) \eqs,
\end{split}
\eeq
for all $t \geq 1$ and $i \in [N]$.
\end{Lemma}
The proof can be found in Appendix~\ref{pf:prob}.
Note that the above is given in terms of $\xi_t^{-1} \bargrd{t}{i}$
instead of $\bargrd{t}{i}$. However, the result remains relevant as
the LO \eqref{eq:lo} in the DeFW algorithm is
\emph{scale invariant}, i.e., $\arg \min_{\atom \in \Cset} \langle \bargrd{t}{i}, \atom \rangle =
\arg \min_{\atom \in \Cset} \langle \alpha \bargrd{t}{i}, \atom \rangle$
for any $\alpha > 0$.
In other words, performing the FW step with
$\bargrd{t}{i}$ is equivalent to doing so with $\xi_t^{-1} \bargrd{t}{i}$.
As $\xi_t^{-1} \bargrd{t}{i}$ is an ${\cal O}(1/t)$ approximation to
$N^{-1} \sum_{j=1}^N \grd f_j ( \bar{\prm}_t^j )$,
Assumption~\ref{ass:a2} is satisfied with $\dgrd_t = {\cal O} (1/t)$.
Lastly, we conclude that
\begin{Corollary}
The sparsified DeFW algorithm using {random} coordinate selection,
i.e., with line~\ref{line:con} \& \ref{line:agg} 
in Algorithm~\ref{alg:defw} replaced by \eqref{eq:gac2} \& \eqref{eq:gac3}, respectively,
generates iterates that satisfy the guarantees in Theorem~\ref{thm:cvx} (with high probability).
Under strong convexity and interior optimal point assumption,
the communication complexity is
${\cal O} ( N \cdot (1/\delta) \cdot \log (1 / \delta ) )$ to reach a $\delta$-optimal solution to \eqref{eq:lasso}.
\end{Corollary}
In the above, the first statement is a consequence of Lemma~\ref{lem:prob}.
The second statement can be verified by noting that
reaching a $\delta$-optimal solution requires ${\cal O} ( 1/ \sqrt{\delta})$ iterations
and the communication cost
is ${\cal O} ( N t \log t )$ at iteration $t$, as
the agents exchange an ${\cal O} ( Nt  )$-sparse vector for 
$\Theta( \log t )$ times.

\section{Numerical Experiments} \label{sec:num}
We perform numerical experiments to verify our theoretical findings on the DeFW algorithm.
The following discussions will focus on the two applications described in Section~\ref{sec:ext}
using synthetic and real data.
To simulate the decentralized optimization setting, we artificially construct
a network of $N=50$ agents, where
the underlying communication network $G$ is an Erdos-Renyi graph with
connectivity of $0.1$. For the AC steps \eqref{eq:gac2}, \eqref{eq:gac1} \& \eqref{eq:gac3},
the doubly stochastic matrix ${\bm W}$
is calculated according to the Metropolis-Hastings rule in \cite{Xiao2004}.

\subsection{Decentralized Matrix Completion}
This section considers the decentralized matrix completion problem,
where the goal is to predict missing entries of an unknown matrix 
through corrupted
partial measurements.
We consider two datasets --- the first dataset is
synthetically generated where the 
unknown matrix $\prm_{\sf true}$ is rank-$K$ and
has dimensions of $m_1 \times m_2$; the matrix is generated as
${\prm}_{\sf true} = \sum_{i=1}^{K} {\bm y}_i {\bm x}_i^\top / K$
where ${\bm y}_i, {\bm x}_i$ have i.i.d.~${\cal N}(0,1)$ entries and different settings of $m_1,m_2,K$
will be experimented.
The second dataset is the \texttt{movielens100k} 
dataset \cite{movielens15}. The unknown matrix $\prm_{\sf true}$
records the movie ratings of $m_1 = 943$ users on $m_2 = 1682$ movies; and a
total of $10^5$ entries in $\prm_{\sf true}$ are available as integers ranging from
$1$ to $5$.
The datasets are divided into training and testing
sets and the mean square error (MSE) on the \emph{testing set} is evaluated as:
\beq \textstyle
{\rm MSE} = | \Omega_{\sf {test}} |^{-1} \sum_{ (k,l) \in \Omega_{\sf {test}} } \big| [\prm_{\sf true}]_{k,l} - [ \hat{\prm}]_{k,l} \big|^2 \eqs,
\eeq
where $ \hat{\prm}$ denotes the estimated $\prm$ produced by the algorithm.

\begin{figure}[t]
\centering
\ifconfver
\includegraphics[width=.5\linewidth]{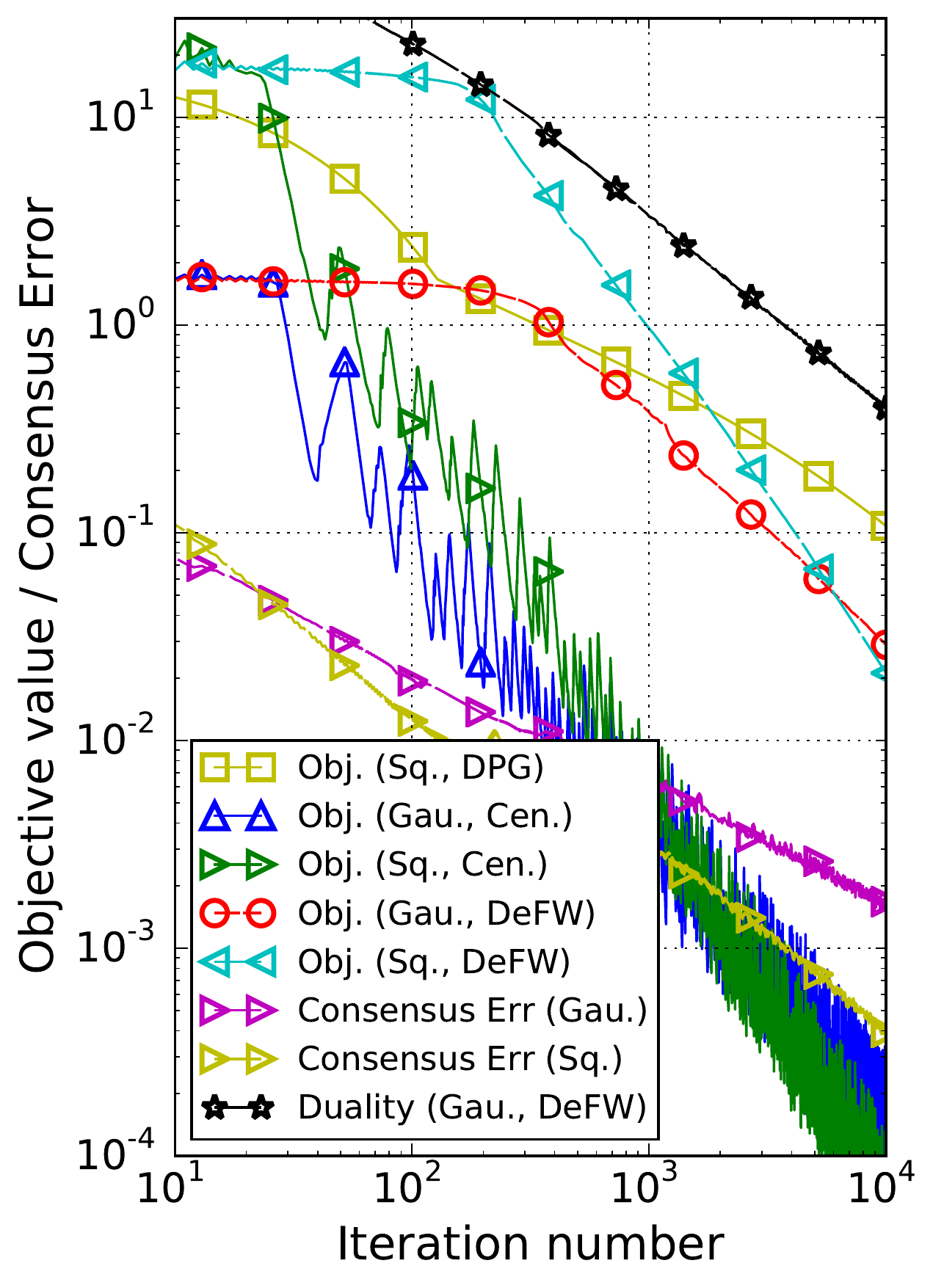}\includegraphics[width=.5\linewidth]{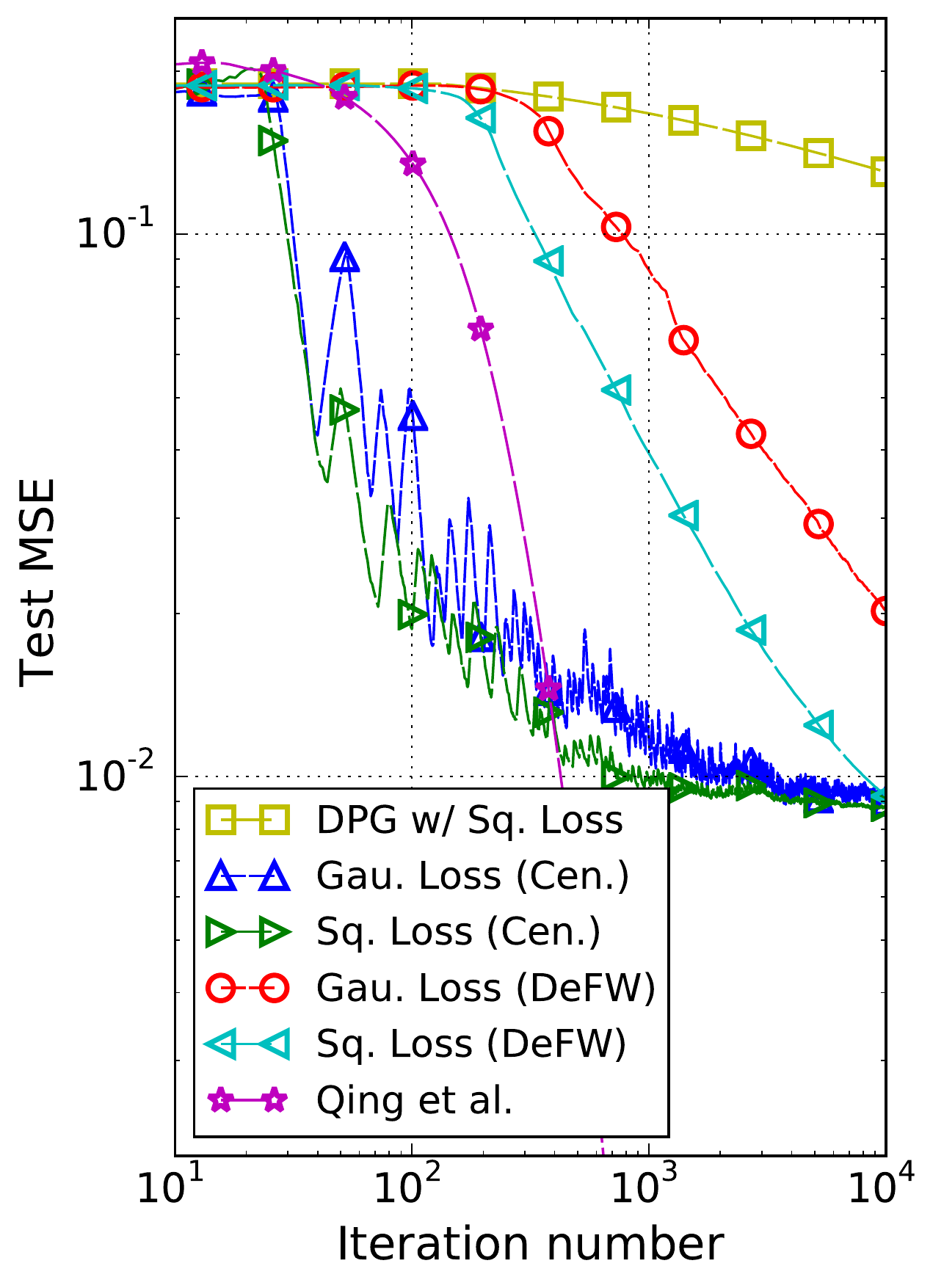} \vspace{-.3cm}

\includegraphics[width=.75\linewidth]{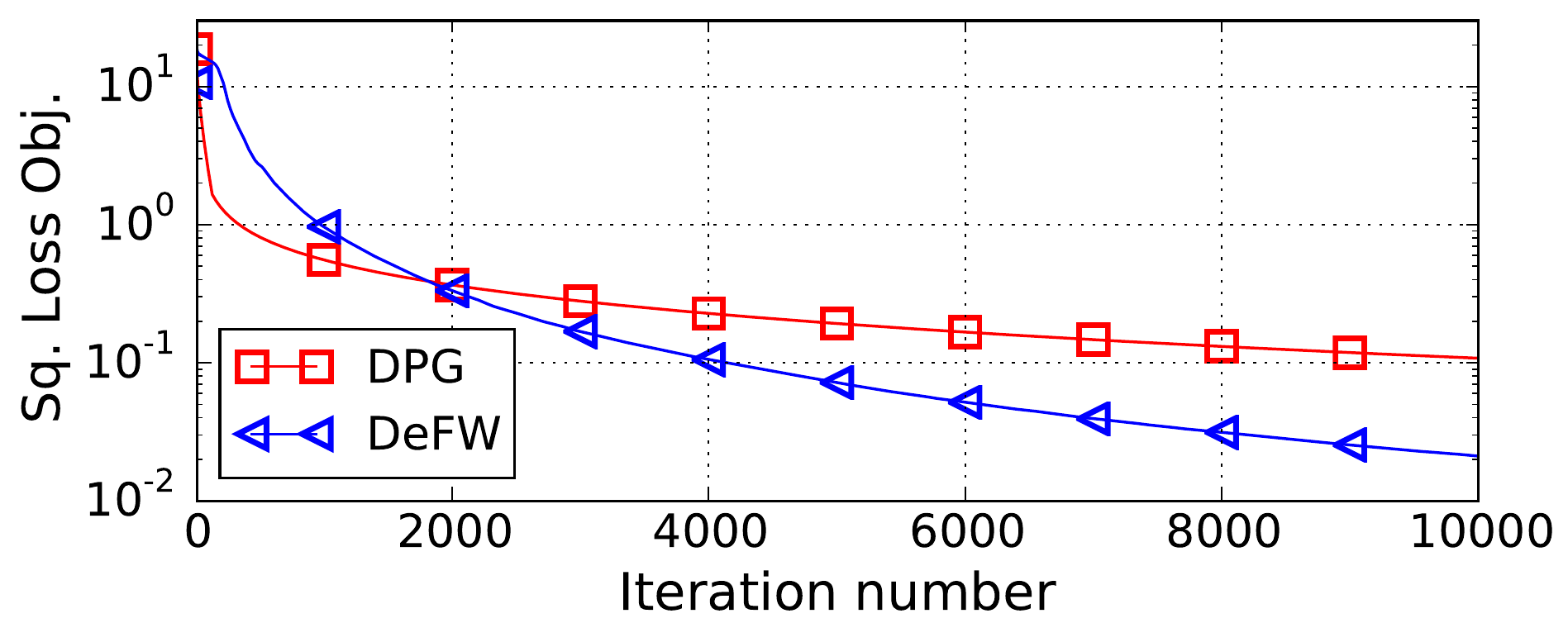}
\vspace{-.3cm}
\else
\includegraphics[width=.36\linewidth]{J_ObjSynData_NoOutlier_DPG_Step.pdf}~\includegraphics[width=.36\linewidth]{J_MSESynData_NoOutlier_DPG_Step.pdf} \vspace{-.2cm}
\includegraphics[width=.6\linewidth]{J_ObjSynData_JustDPG_Step.pdf}
\fi
\caption{Results on noiseless synthetic data with $m_1=100, m_2=250$ and rank $K=5$. (Top-Left) Objective value and consensus error of $\bar{\prm}_t^i$ against iteration number $t$, the objective values are evaluated by $F(\bar{\prm}_t)$. (Top-Right) Worst-case MSE (among agents) against iteration number on the testing set. (Bottom) Objective value (sq.~loss)
against iteration number $t$ for
DeFW and DPG. The legend `Gau.', `Sq.' denote the consensus-based DeFW algorithm applied to \eqref{eq:mc}
with the negated Gaussian and square loss, respectively.} \vspace{-.2cm}
\label{fig:syn_lowrank}
\end{figure}

For the synthetic dataset, the training (testing) set contains $20\%$ ($80\%$) entries
which are selected randomly. For
\texttt{movielens100k}, the training (testing) set contains $80 \times 10^3$ ($20 \times 10^3$) entries.
The training data of the two datasets are equally partitioned into $N=50$ parts; 
for \texttt{movielens100k}, each
agent holds $1600$ entries.
We evaluate the performance of the proposed consensus-based DeFW algorithm
applied to 
square loss and negated Gaussian loss, as described in Section~\ref{sec:mc}.
Unless otherwise specified, we fix the number of AC rounds applied at $\ell=1$ such
that the agents only exchange information once per iteration.
As the negated Gaussian loss is non-convex, we set the step size as $\gamma_t = t^{-0.75}$.
The centralized FW algorithm for both losses will also be compared (cf.~\eqref{eq:fw}); as well as the decentralized
algorithm in \cite{qing12} (labeled as `Qing et al.')
and the DPG algorithm \cite{RNV2010} with step size set to $\alpha_t = 0.1 N / ( \sqrt{t}+1)$ applied to square loss.

Our first example considers the noiseless synthetic dataset of problem dimension $m_1=100, m_2=250$, $K=5$.
The results are shown in Fig.~\ref{fig:syn_lowrank}.
Here, for the DeFW/DPG algorithms, we set the trace-norm radius to $R = 1.2 \| \prm_{\sf true} \|_{\sigma,1}$;
and the algorithm in \cite{qing12} is supplied with the true rank $K$ of $\prm_{\sf true}$.
Notice that for this set of data, the minimum of \eqref{eq:mc} can be achieved by
$\prm = \prm_{\sf true} \in \Cset$ with a zero objective value.
From the top-left plot, for the DeFW algorithm applied to the convex square loss function,
we observe an ${\cal O}(1/t^2)$ trend for the objective values, corroborating
with our analysis in Theorem~\ref{thm:cvx}; for the non-convex  negated Gaussian loss function,
the objective value and the FW/duality gap $g_t$
also decay with $t$, the latter indicates the convergence to a stationary point.
Moreover, the consensus error of $\bar{\prm}_t^i$ for DeFW applied to the two objective functions
decay at the rate predicted by Lemma~\ref{lem:gac2}.
On the other hand, the top-right plot compares mean square error (MSE) of the predicted matrix $\prm$
for the testing set.
Here, we also compare the result with the algorithm in \cite{qing12}.
We observe that the MSE performance of the DeFW algorithms approach their centralized counterpart
as the iteration number grows, yet the algorithm in \cite{qing12} achieves the best performance in
this setting, notice that the true rank of $\prm_{\sf true}$ is provided to this algorithm.
From the bottom plot, the DPG method applied to square loss function converges 
at a relatively fast rate in the beginning, 
but was overtaken by DeFW as the iteration number grows. 
It is worth mentioning that the DeFW algorithms have a consistently
better  MSE performance than DPG.

\begin{figure}[t]
\centering
\ifconfver
\includegraphics[width=.525\linewidth]{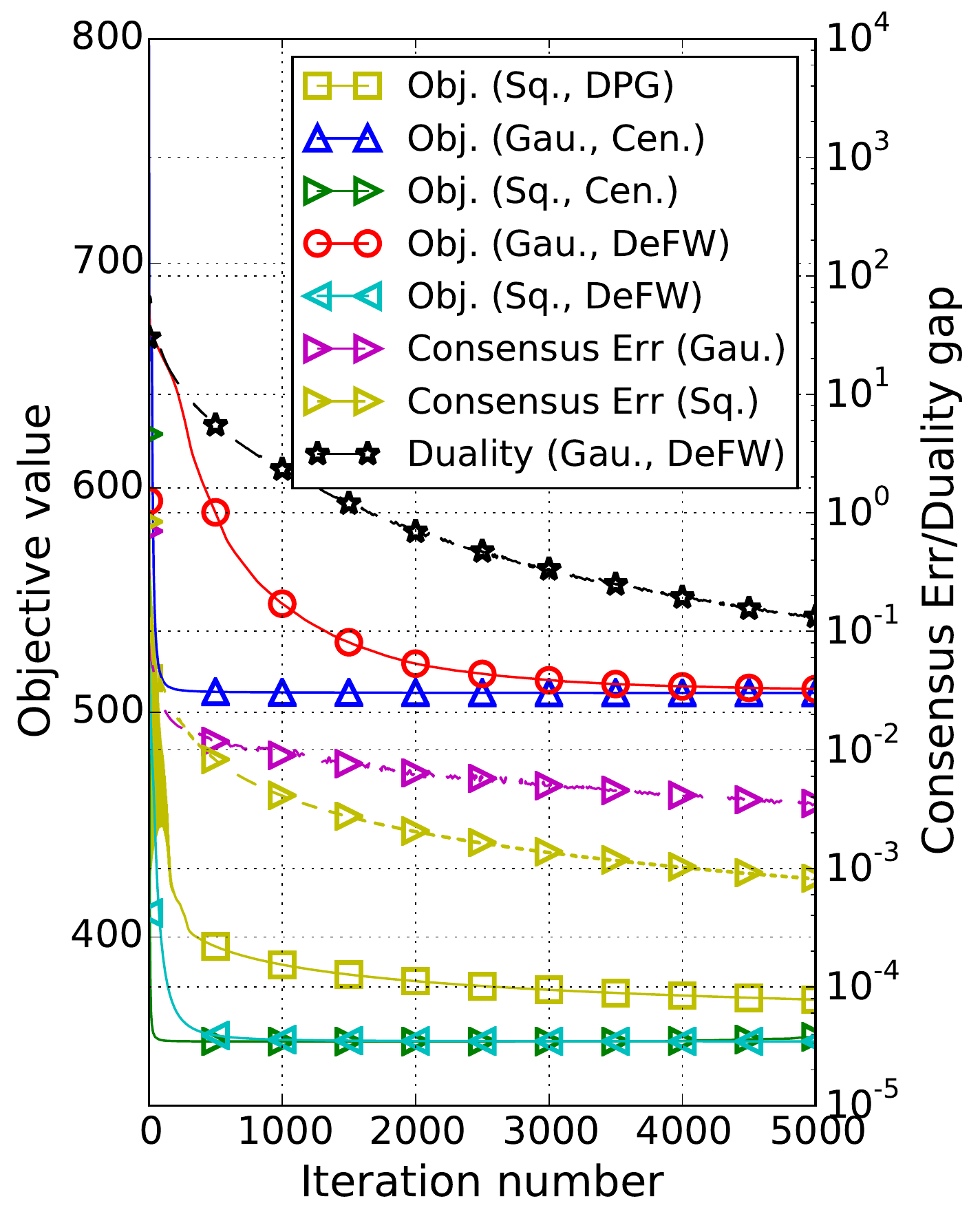}\includegraphics[width=.475\linewidth]{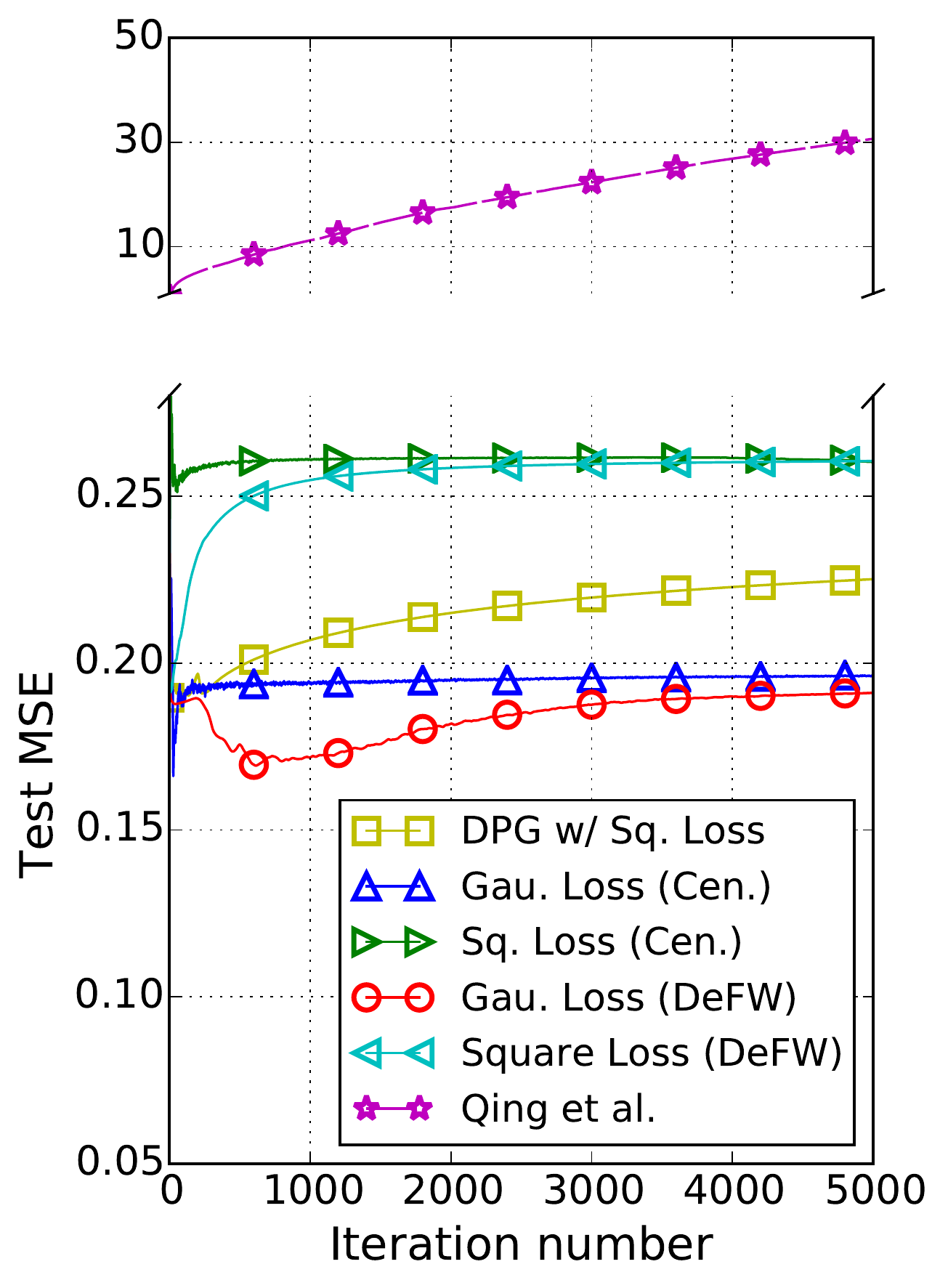} \vspace{-.5cm}
\else
\includegraphics[width=.38\linewidth]{J_ObjSynData_YesOutlier_Rev_Step.pdf}~\includegraphics[width=.35\linewidth]{J_MSESynData_YesOutlier_Rev_Step.pdf} \vspace{-.3cm}
\fi
\caption{Results on sparse-noise contaminated synthetic data with $m_1=100, m_2=250$ and rank $K=5$. (Left) Objective value and consensus error of $\bar{\prm}_t^i$ against the DeFW iteration number $t$. Notice that the consensus error (in purple and yellow) / duality gap (in black) are plotted in a logarithmic scale (cf.~the right y-axis) while the objective values are plotted in a linear scale; (Right) MSE against the DeFW iteration number $t$ on the testing set.}
\vspace{-.3cm}\label{fig:syn_outlier}
\end{figure}

The second example considers adding noise to the observations
for the same synthetic data case in Fig.~\ref{fig:syn_lowrank}. 
In particular, we adopt the
same setting as the previous example but include a \emph{sparse}
noise in the observations --- here, each $Z_{k,l} = p_{k,l} \cdot \tilde{Z}_{k,l}$
where $p_{k,l}$ is Bernoulli with $P( p_{k,l} = 1 ) = 0.2$ and
$\tilde{Z}_{k,l} \sim {\cal N}(0,5)$ (cf.~\eqref{eq:mat_obs}).
The convergence results are compared in Fig.~\ref{fig:syn_outlier}.
For the left plot, we observe similar convergence behaviors for the DeFW algorithms
applied to different objective functions as in the previous example.
On the right plot, we observe that the DeFW algorithm based on negated Gaussian loss achieves the lowest
MSE, demonstrating its robustness to outlier noise.
We also see that the algorithm in \cite{qing12} performs poorly on this dataset.

\begin{figure}[t]
\centering
\ifconfver
\includegraphics[width=.7\linewidth]{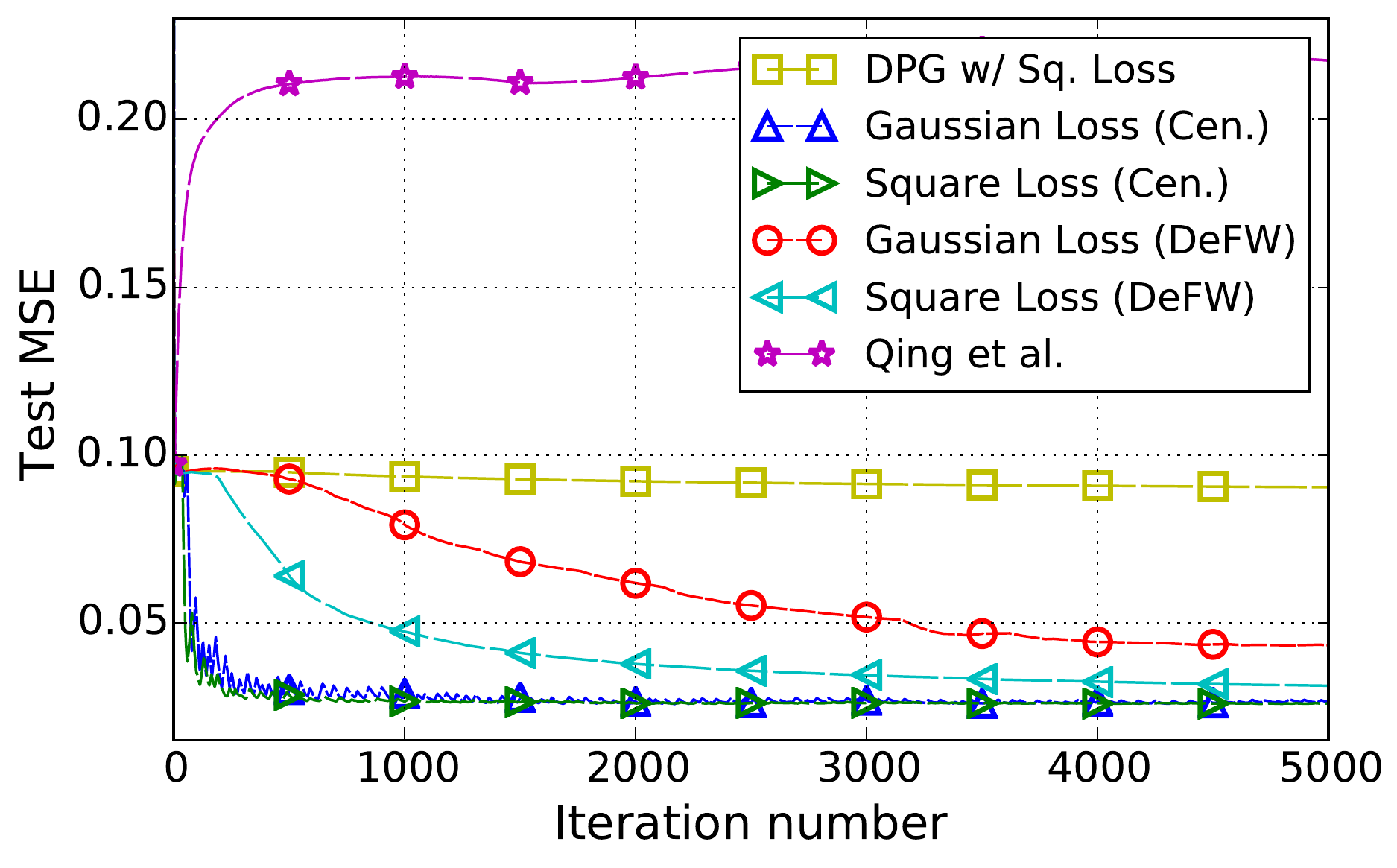} \vspace{-.3cm}
\else
\includegraphics[width=.5\linewidth]{J_MSESynData_R10_Rev_Step.pdf} \vspace{-.3cm}
\fi
\caption{Convergence of test MSE against iteration number on the testing set on noise-free synthetic data with $m_1=100, m_2=250$ and rank $K=10$.}\vspace{-.2cm}
\label{fig:syn_R10}
\end{figure}

Another interesting discovery is that the algorithm in \cite{qing12} seems to
fail when the rank of $\prm_{\sf true}$
is high, even when the true rank is known and the observations are noiseless.
In Fig.~\ref{fig:syn_R10}, we show the MSE against iteration number
of the algorithms when the synthetic data is noiseless and generated with $m_1=100,m_2=250,K=10$.
As seen, \cite{qing12} fails to produce a low MSE,
while DeFW offers a reasonable performance.

\begin{figure}[t]
\centering
\ifconfver
\includegraphics[width=.5\linewidth]{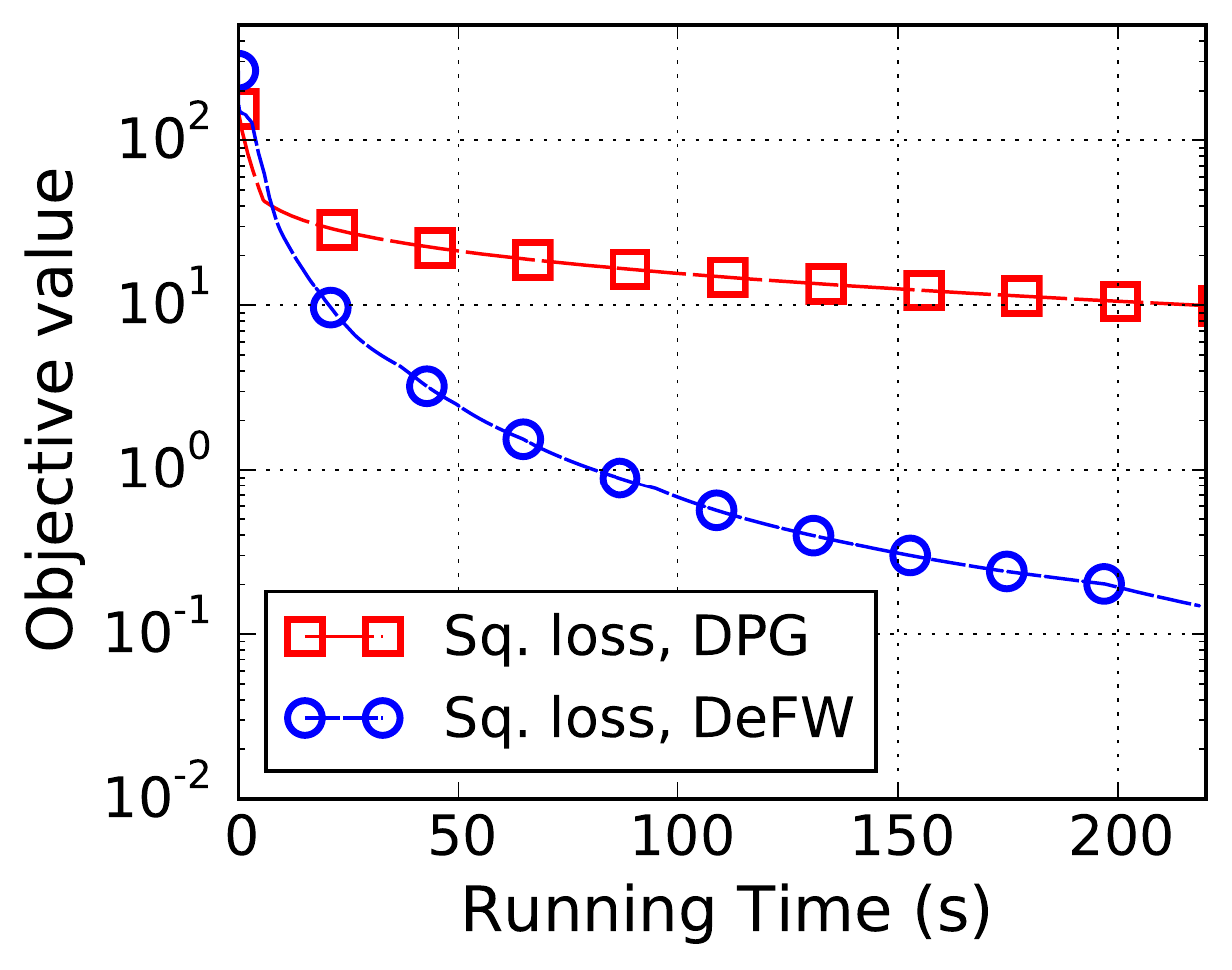}\includegraphics[width=.5\linewidth]{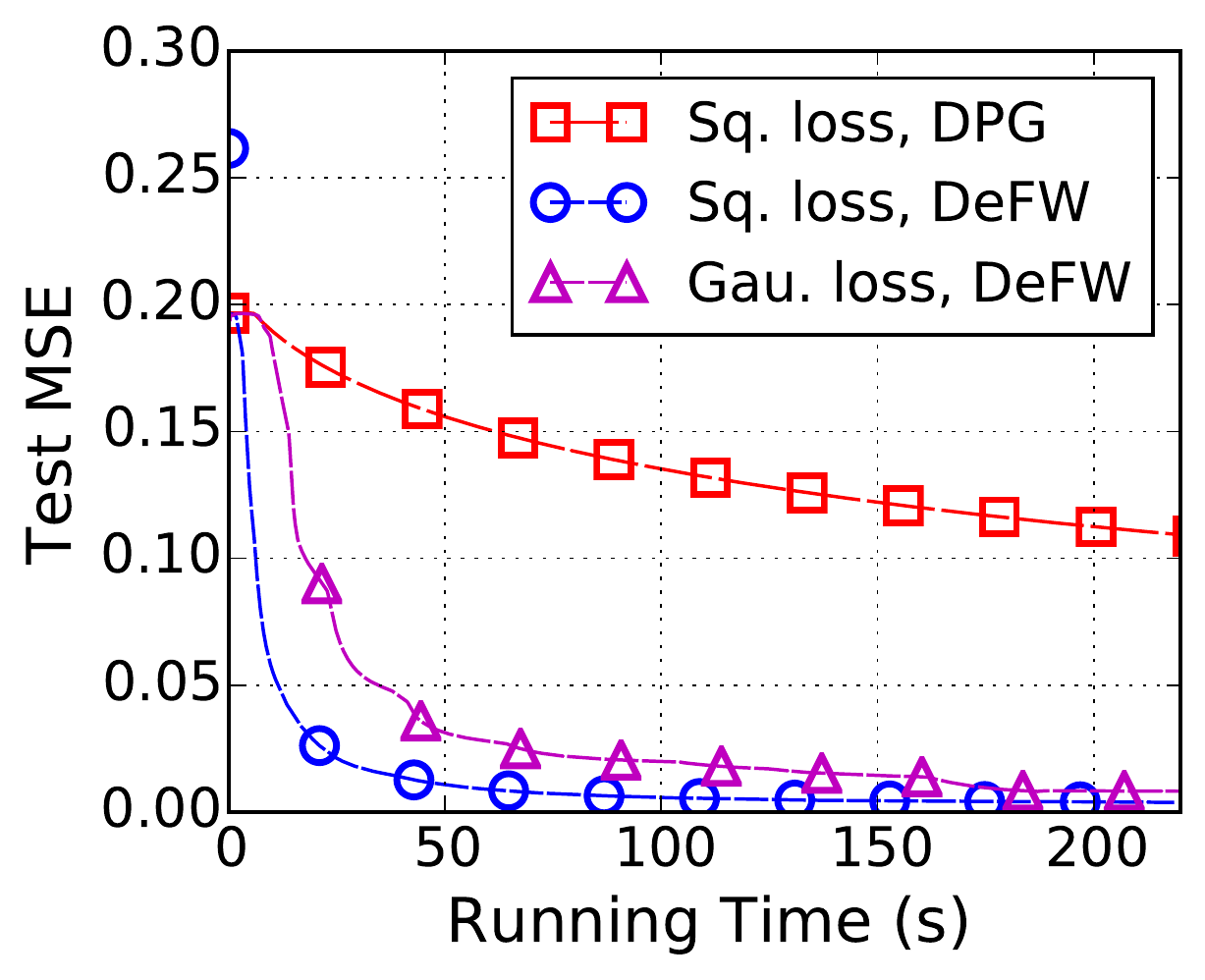} \vspace{-.6cm}
\else
\includegraphics[width=.38\linewidth]{J_TimeObj_DPG_Step.pdf}\includegraphics[width=.38\linewidth]{J_TimeMSE_DPG_Step.pdf} \vspace{-.3cm}
\fi
\caption{  Results on noiseless synthetic data with $m_1=200, m_2=1000$ and rank $K=5$. (Left) Objective value against running time. (Right) Worst-case MSE (among agents) against running time.} \vspace{-.2cm}
\label{fig:syn_time}
\end{figure}

The next example evaluates the objective value and test MSE on synthetic, noiseless data
against the average runtime per agent. 
We focus on comparing the DeFW and DPG algorithms.  
In Fig.~\ref{fig:syn_time}, DeFW 
demonstrates a significant advantage over DPG 
since the former does not require the projection computation.
In fact, the average running time \emph{per iteration} 
of DeFW is five times faster than 
DPG. 
We also expect the complexity advantages to widen as the problem size grows.

\begin{figure}[t!]
\centering
\ifconfver
\includegraphics[width=.49\linewidth]{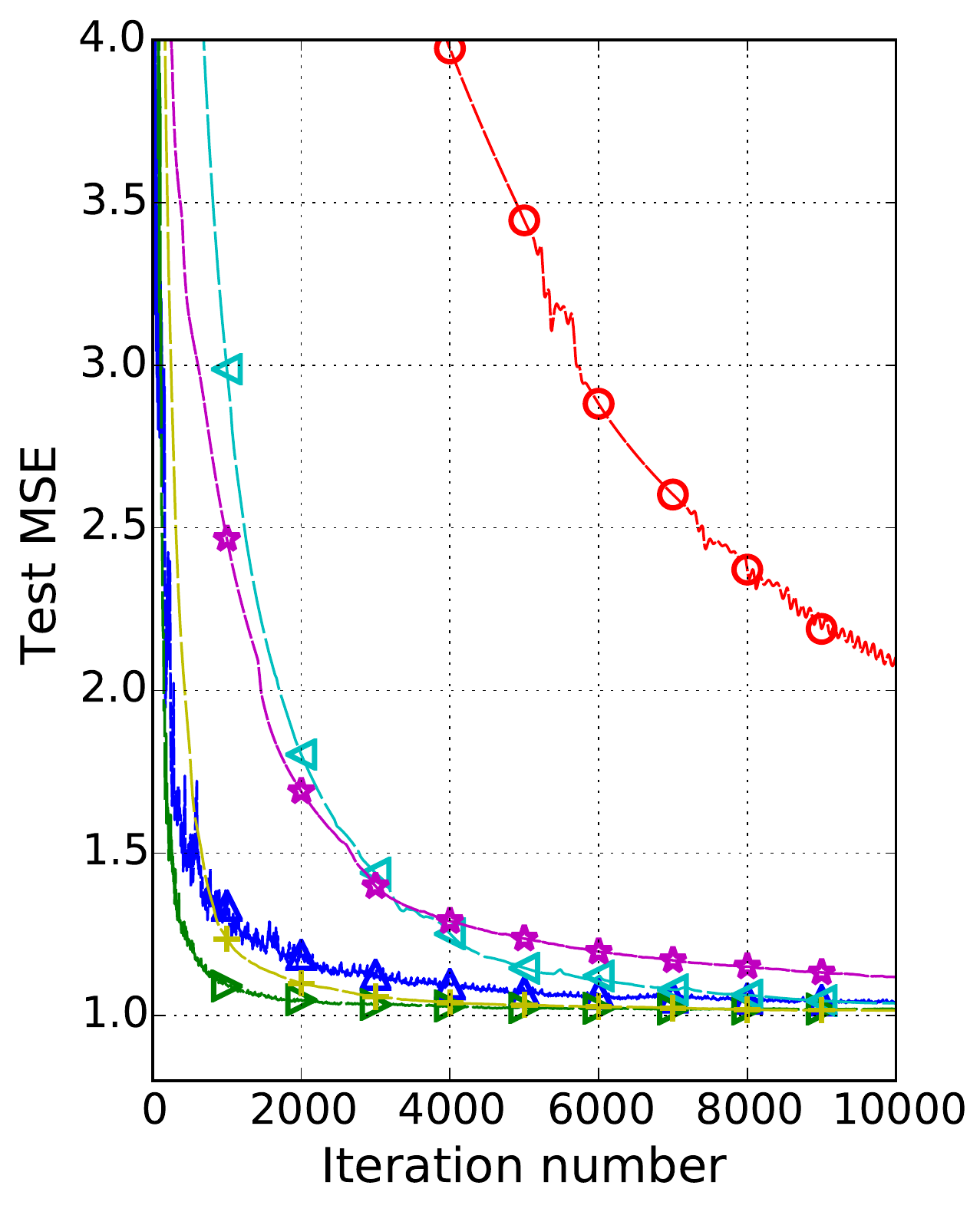}\includegraphics[width=.505\linewidth]{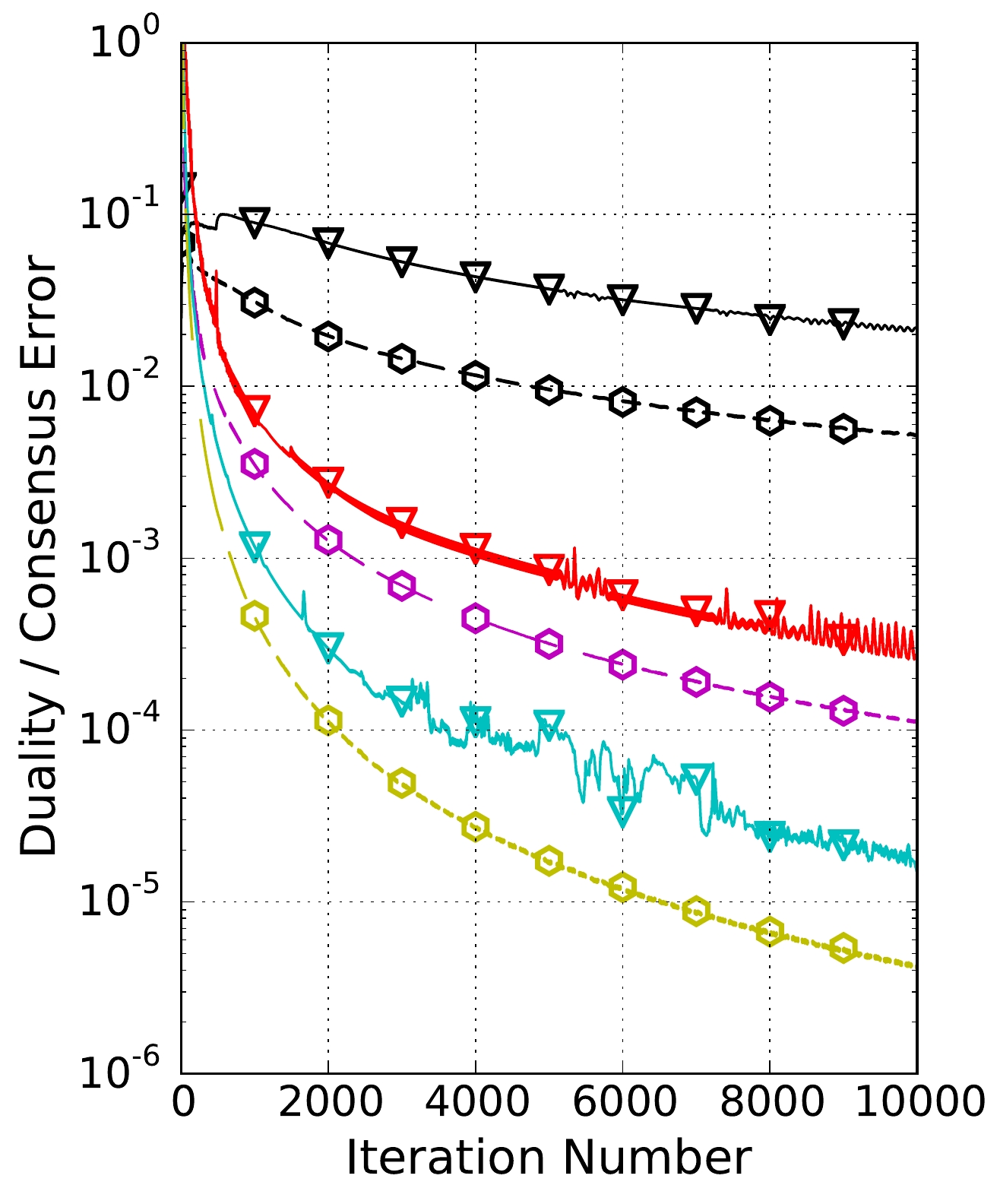} \vspace{-.3cm} \\
\includegraphics[width=.49\linewidth]{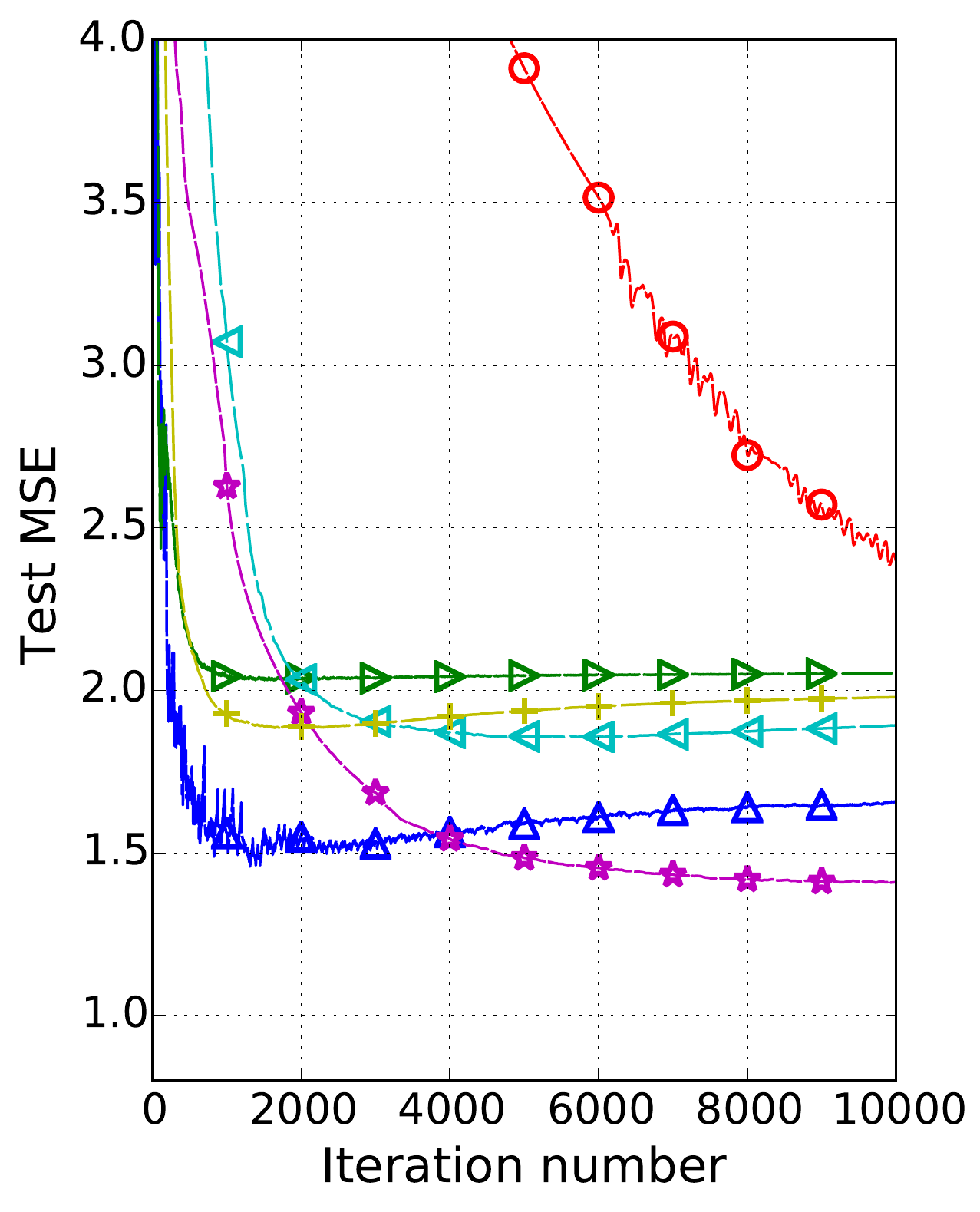}\includegraphics[width=.505\linewidth]{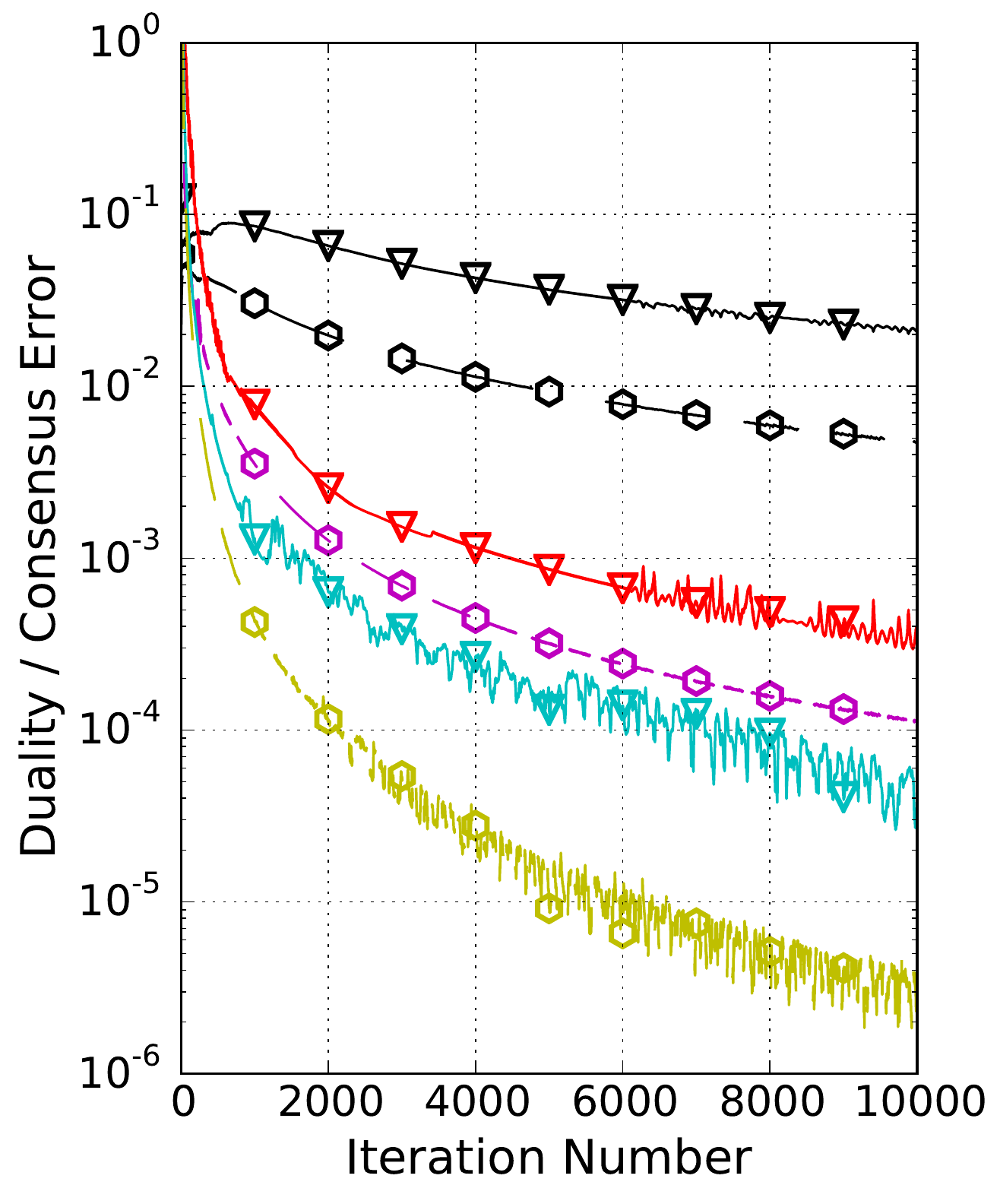} \vspace{-.2cm} \\
\includegraphics[width=.99\linewidth]{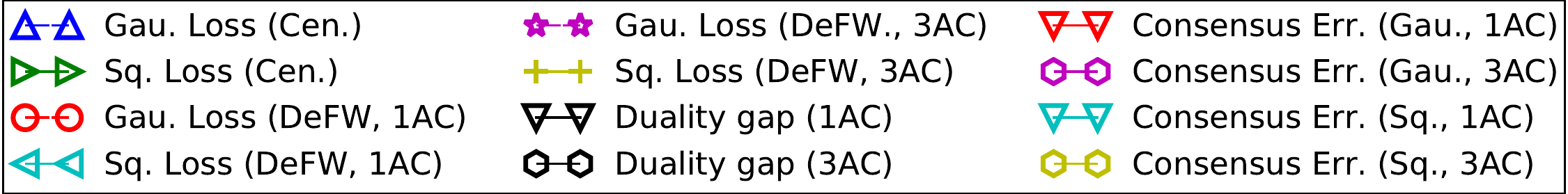}
\else
\includegraphics[width=.35\linewidth]{J_MSEML_NoOutlier_MSEOnly_Rev.pdf}\includegraphics[width=.355\linewidth]{J_MSEML_NoOutlier_DualityOnly_Rev.pdf} \vspace{-.0cm} \\
\includegraphics[width=.35\linewidth]{J_MSEML_YesOutlier_MSEOnly_Rev.pdf}\includegraphics[width=.355\linewidth]{J_MSEML_YesOutlier_DualityOnly_Rev.pdf} \vspace{-.0cm} \\
\includegraphics[width=.7\linewidth]{Legend_ML100k.pdf}
\fi
\caption{Convergence of the DeFW algorithm on \texttt{movielens100k}. (Top) Noiseless observations; (Bottom) sparse-noise contaminated observations. Note that the duality gap and consensus error are drawn in a logarithmic scale in the right plots.}\vspace{-.4cm}
\label{fig:ML}
\end{figure}

Lastly, we consider the dataset \texttt{movielens100k}.
We set $R= 10^5$ and focus on the test MSE evaluated 
against the iteration number
for the proposed DeFW algorithm.
The numerical results are presented in Fig.~\ref{fig:ML}, where we also
compare the case when we apply multiple ($\ell =1, 3$) rounds of AC updates
per iteration to speed up the algorithm.
The left plot in Fig.~\ref{fig:ML} considers the noiseless scenario.
As seen, the proposed DeFW algorithm applied on different loss functions converge
to a reasonable MSE that is attained by the centralized FW algorithm.
We see that the DeFW with negated Gaussian loss has a slower convergence compared to the square
loss which is possibly attributed to the non-convexity of the loss function.
Moreover, the algorithms achieve much faster convergence if we allow $\ell=3$ AC rounds
of  network information exchange per iteration.
The right plot in Fig.~\ref{fig:ML} considers when the observations are
contaminated with a sparse noise
of the same model as Fig.~\ref{fig:syn_outlier}.
We observe that the negated Gaussian loss implementations attain the best MSE
as the non-convex loss is more robust against the sparse noise.
Interestingly, the DeFW algorithm with $\ell=3$ AC rounds has even outperformed its
centralized counterpart.
We suspect that this is caused by the DeFW algorithm converging to a different
local minima for the non-convex problem.

\subsection{Communication-efficient LASSO}
This section conducts numerical experiments on the decentralized sparse learning problem.
We focus on the \emph{sparsified DeFW algorithm} in
Section~\ref{sec:lasso} that has better communication efficiency.
We evaluate the performance  
on both synthetic and benchmark data.
For the synthetic data,
we randomly generate each ${\bm A}_i$ as a $(m=20) \times (d=10000)$ matrix with
${\cal N}(0,1)$ elements (cf.~\eqref{eq:sparse}) and $\prm_{\sf true}$ is a random sparse vector with $\| \prm_{\sf true} \|_0 = 50$
such that the non-zero elements are also ${\cal N}(0,1)$.
The observation noise ${\bm z}_i$ has a variance of $\sigma^2 = {0.01}$.
For benchmark data, we test our method on \texttt{sparco7}
\cite{sparco07}, which is a commonly used
dataset for benchmarking sparse recovery algorithms.
For \texttt{sparco7}, we have ${\bm A}_i \in \RR^{12 \times 2560}$
as the local measurement matrix and $\prm_{\sf true}$ is a sparse vector with $\| \prm_{\sf true} \|_0 = 20$.

The sparsified DeFW algorithm is implemented with $p_t = \lceil 2 + \alpha_{comm} \cdot t \rceil$,
 $\ell_t = \lceil \log(t) + 1 \rceil$ with extreme or random coordinate selection.
We compare the algorithms of
PG-EXTRA \cite{Shi2015} (with fixed step size $\alpha = 1/d$), 
DPG \cite{RNV2010} (with step size $\alpha_t = 1 / {t}$)  and BHT \cite{ravazzi15}.
 DeFW, PG-EXTRA and DPG are set to solve the convex problem \eqref{eq:lasso}
with $R = 1.1 \| \prm_{\sf true} \|_1$.
 BHT is a communication efficient decentralized version of
IHT and is supplied with the true sparsity level in
our simulations.

\begin{figure}[t!]
\centering
\ifconfver
\includegraphics[width=.5\linewidth]{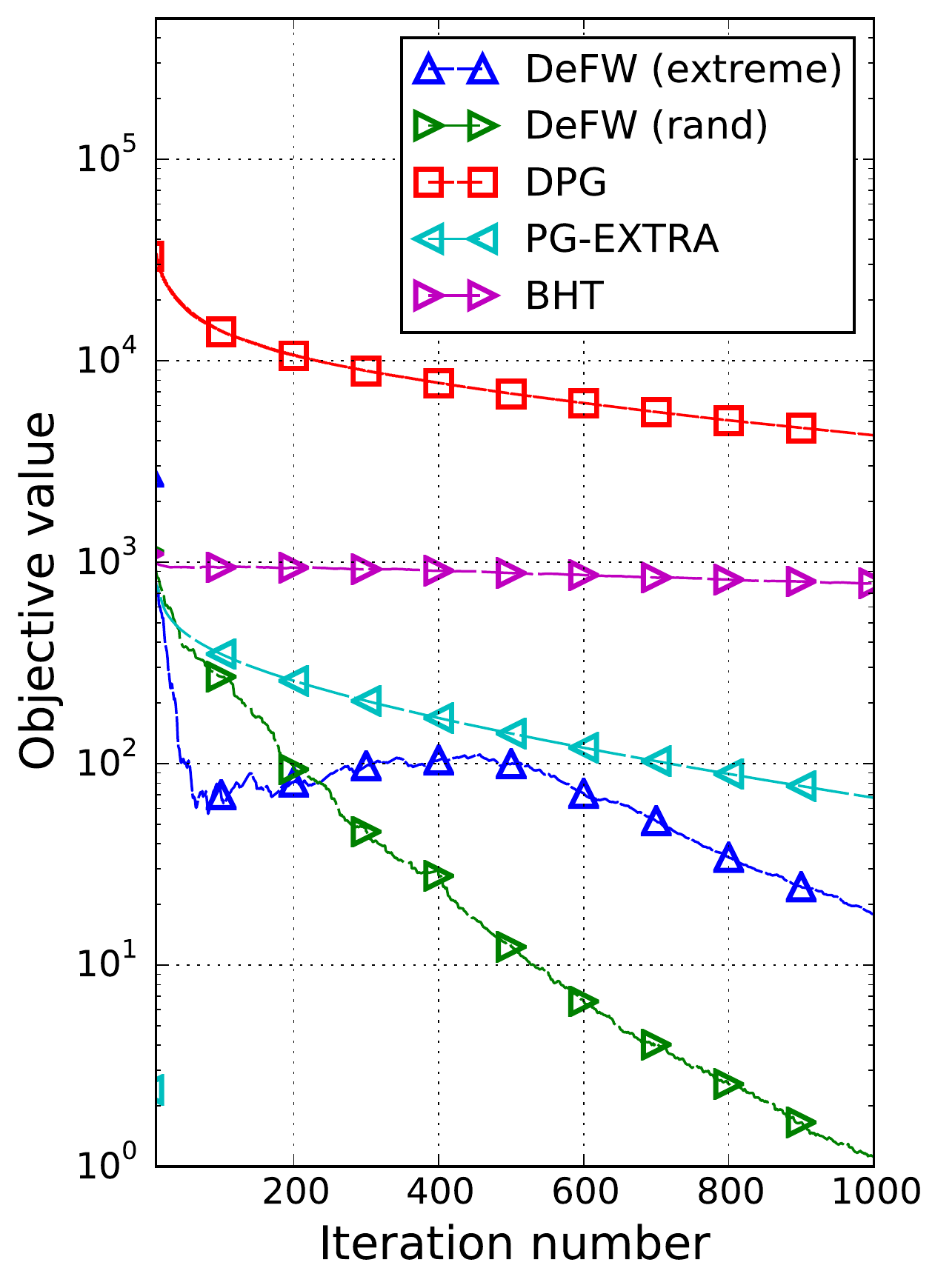}\includegraphics[width=.475\linewidth]{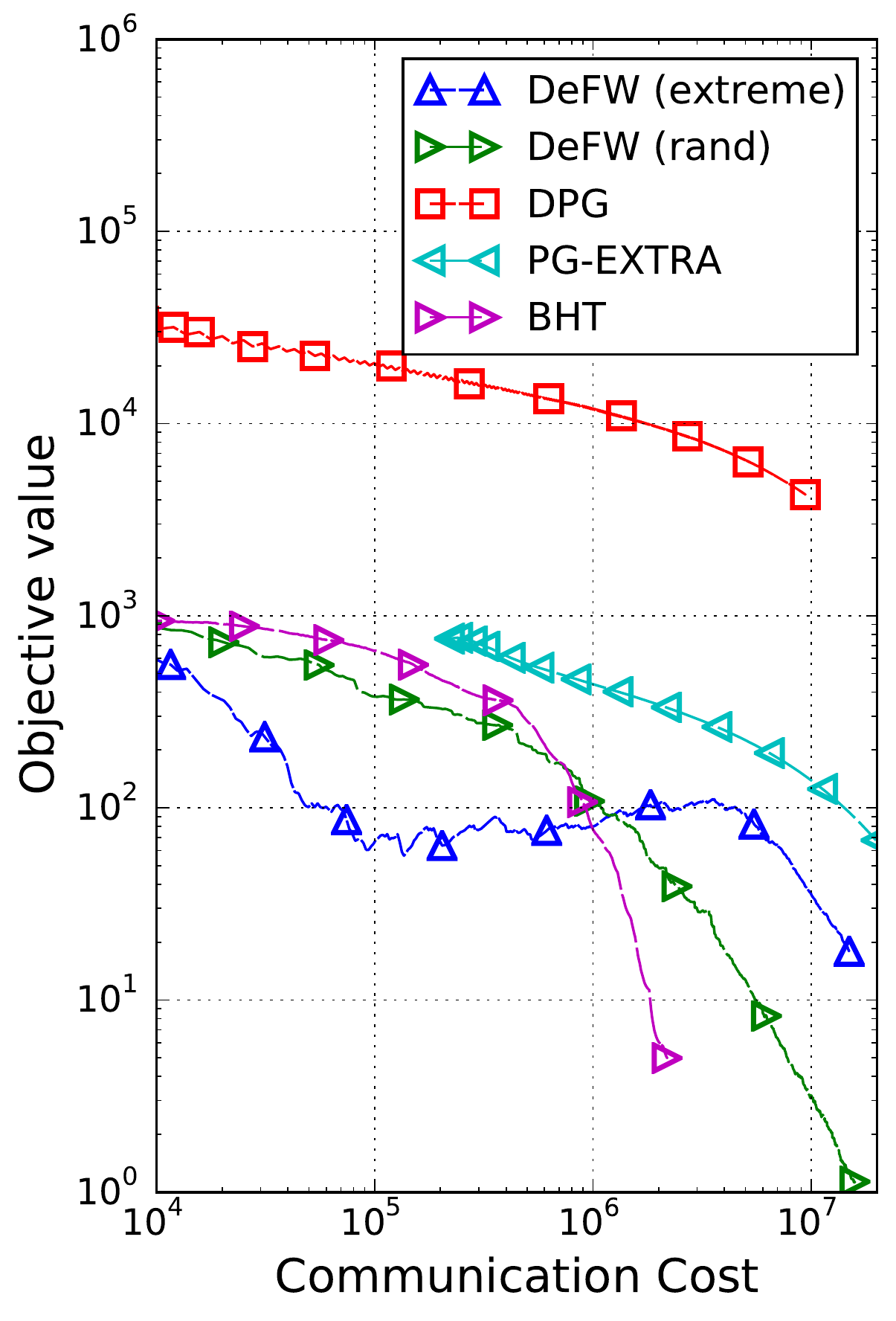} \vspace{-.3cm}
\else
\includegraphics[width=.35\linewidth]{J_LASSO_Syn_Rev.pdf}~\includegraphics[width=.33\linewidth]{J_LASSO_Comm_Syn_Rev.pdf} \vspace{-.3cm}
\fi
\caption{Convergence of the objective value on LASSO with synthetic dataset. (Left) against the iteration number. (Right) against the communication cost (i.e., total number of values transmitted/received in the network during AC updates). In the legend, `DeFW (extreme)' refers to the extreme coordinate selection and `DeFW (rand)' refers to the random coordinate selection scheme.}
\label{fig:lasso}
\end{figure}

The first example in Fig.~\ref{fig:lasso} shows the convergence of the algorithms
on the synthetic data, where
we compare the objective value against the number of iterations and the communication cost,
i.e., total number of values sent during the distributed optimization.
We set $\alpha_{comm} = 0.05$ for the DeFW algorithms.
From the left plot, we observe that DeFW and PG-EXTRA algorithms have similar
iteration complexity while `DeFW (rand)' seems to have the fastest convergence.
Meanwhile, BHT demands a high number of iterations for convergence.
On the other hand, in the right plot, the DeFW algorithms
demonstrate the best communication efficiency at low accuracy,
while they lose to BHT at higher accuracy.
Lastly, `DeFW (extreme)' achieves a better accuracy at the beginning 
(i.e., less communication cost paid) but is overtaken by `DeFW (rand)' as the communication
cost grows.

\begin{figure}[t]
\centering
\ifconfver
\includegraphics[width=.8\linewidth]{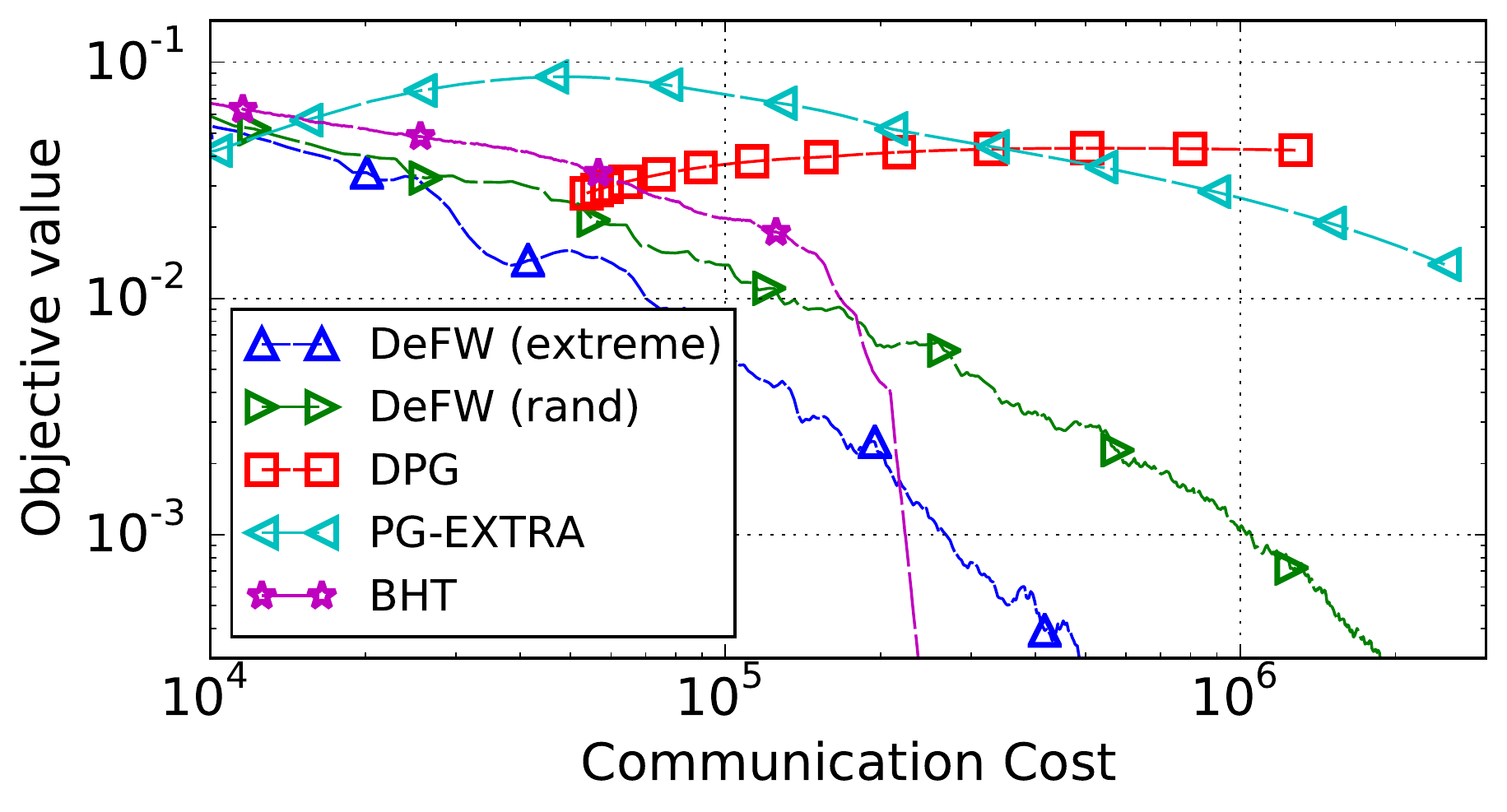} \vspace{-.3cm}
\else
\includegraphics[width=.5\linewidth]{J_LASSO_Comm_Sparco7_Rev.pdf} \vspace{-.3cm}
\fi
\caption{Convergence of the objective value against the communication cost on LASSO with \texttt{sparco7} dataset. In the legend, `DeFW (extreme)' refers to the extreme coordinate selection and `DeFW (rand)' refers to the random coordinate selection scheme.}\vspace{-.2cm}
\label{fig:sparco7}
\end{figure}

We then compare the performance on \texttt{sparco7}, where we show the
convergence of objective value against the communication cost in Fig.~\ref{fig:sparco7}.
We set $\alpha_{comm} = 0.025$ for the sparsified DeFW algorithms.
At low accuracy, the  DeFW algorithms
offer the best communication cost-accuracy trade-off, i.e., it performs the best at
an accuracy of above $\sim 10^{-2}$.
Moreover, `DeFW (extreme)' seems to be perform better than `DeFW (rand)'
in this example.
Nevertheless, the BHT algorithm achieves the best performance when
the communication cost paid is above $3 \times 10^5$.
Lastly, we comment that although BHT has
the lowest communication cost at \emph{high} accuracy, its computational
complexity is high as the former requires a large number of iterations
to reach a reasonable accuracy (cf.~left plot of Fig.~\ref{fig:lasso}).
The sparsified DeFW offers a better balance of the
communication and computation complexity.

\section{Conclusions \& Open Problems}
In this paper, we have studied a decentralized projection-free algorithm for constrained
optimization, which we called the DeFW algorithm. Importantly, we showed that the
DeFW algorithm converges for both convex and non-convex loss functions and the
respective convergence rates are analyzed.
The efficacy of the proposed algorithm is demonstrated through tackling two problems related
to machine learning, with the advantages over previous state-of-the-art
demonstrated through numerical experiments.
Future directions of study include developing an asynchronous version of the DeFW algorithm.

\section{Acknowledgement}
The authors would like to thank the anonymous reviewers
for providing constructive comments on our paper.

\ifplainver
\appendix
\else
\appendices
\fi

\section{Proof of Theorem~\ref{thm:cvx}} \label{pf:cvx}
The proof of Theorem~\ref{thm:cvx} follows from our recent analysis on  online/stochastic
FW algorithm \cite{ourpaper}.
Using line~\ref{line:fw} of Algorithm~\ref{alg:defw}, we observe that:
\beq
\textstyle
\bar{\prm}_{t+1} = \bar{\prm}_t + \gamma_t ( N^{-1} \sum_{i=1}^N \atom_t^i - \bar{\prm}_t ) \eqs.
\eeq
Define $h_t \eqdef F(\bar{\prm}_t) - F(\prm^\star)$ where $\prm^\star$ is an optimal solution
to \eqref{eq:opt}. From the $L$-smoothness of $F$ and the boundedness of $\Cset$, we have:
\beq \label{eq:fw_key}
h_{t+1} \leq h_t + \frac{\gamma_t}{N} \sum_{i=1}^N \langle \atom_t^i - \bar{\prm}_t , \grd F(\bar{\prm}_t ) \rangle + \gamma_t^2 \frac{L \bar{\rho}^2}{2}\eqs,
\eeq
where $\bar{\rho}$ was defined in \eqref{eq:rho}.
Observe the following chain for the inner product: for each $i \in [N]$, we have
\beq \label{eq:err_bd} \begin{split}
\langle \atom_t^i - \bar{\prm}_t , \grd F(\bar{\prm}_t ) \rangle & \leq \langle \atom_t^i - \bar{\prm}_t , \bargrd{t}{i} \rangle + \bar{\rho} \| \bargrd{t}{i} - \grd F(\bar{\prm}_t) \|_2 \\
& \hspace{-2.2cm} \leq \langle \atom - \bar{\prm}_t , \bargrd{t}{i} \rangle + \bar{\rho} \cdot \| \bargrd{t}{i} - \grd F(\bar{\prm}_t) \|_2,~\forall~ \atom \in \Cset \\
& \hspace{-2.2cm} \leq \langle \atom - \bar{\prm}_t , \grd F(\bar{\prm}_t) \rangle + 2 \bar{\rho} \cdot \| \bargrd{t}{i} - \grd F(\bar{\prm}_t) \|_2,~\forall~ \atom \in \Cset,
\end{split}
\eeq
where we have added and subtracted $\bargrd{t}{i}$ in the first inequality; and used the fact
$\atom_t^i \in \arg \min_{ \atom \in \Cset } \langle \atom, \bargrd{t}{i} \rangle$ in the second inequality.
Recalling that $\bargrd{t}{}= N^{-1} \sum_{i=1}^N \grd f_i( \bar{\prm}_t^i )$, 
\beq \label{eq:thm1_result}
\begin{split}
\| \bargrd{t}{i} - \grd F(\bar{\prm}_t) \|_2 & \\
& \hspace{-1.8cm} \leq \| \bargrd{t}{i} - \bargrd{t}{} \|_2 + \| \bargrd{t}{} - \grd F(\bar{\prm}_t) \|_2 \\
& \hspace{-1.8cm}\textstyle \leq \dgrd_t + N^{-1} \sum_{i=1}^N \| \grd f_i( \bar{\prm}_t^i ) - \grd f_i (\bar{\prm}_t) \|_2 \\
& \hspace{-1.8cm}\textstyle \leq \dgrd_t + L \cdot N^{-1} \sum_{i=1}^N \| \bar{\prm}_t^i - \bar{\prm}_t \|_2 \\
& \hspace{-1.8cm} \leq \dgrd_t + L \cdot \dprm_t \eqs, \vspace{.4cm}
\end{split}  \vspace{-.4cm}
\eeq
where the third inequality is due to the $L$-smoothness of $\{ f_i \}_{i=1}^N$.
Recalling that $\dprm_t = C_p / t$, $\dgrd_t = C_g / t$ and substituting the results above into the inequality \eqref{eq:fw_key} implies:
\beq \label{eq:fw_key2}
h_{t+1} \leq h_t + \gamma_t \langle \bar{\atom}_t - \bar{\prm}_t , \grd F(\bar{\prm}_t) \rangle
+ \gamma_t^2 \frac{L \bar{\rho}^2}{2} + 2 \bar{\rho} \gamma_t \frac{C_g + L C_p}{t} \eqs,
\eeq
where $\bar{\atom}_t \in \Cset$ is the minimizer of the linear optimization
\eqref{eq:lo} using $\grd F(\bar{\prm}_t)$, i.e.,
\beq \label{eq:abar}
\bar{\atom}_t \in \arg \min_{ \atom \in \Cset } \langle \atom, \grd F(\bar{\prm}_t) \rangle \eqs.
\eeq

\noindent \textbf{Case 1}: When $F$ is convex, we observe
\beq
\langle \bar{\atom}_t - \bar{\prm}_t , \grd F(\bar{\prm}_t) \rangle \leq \langle \prm^\star - \bar{\prm}_t , \grd F(\bar{\prm}_t) \rangle \leq - h_t \eqs,
\eeq
where the first inequality is due to the optimality of $\bar{\atom}_t$ and the
last inequality stems from the convexity of $F$. Plugging the
above into \eqref{eq:fw_key2} yields
\beq \label{eq:slow_conv}
h_{t+1} \leq (1-\gamma_t) h_t + \gamma_t^2 \frac{L \bar{\rho}^2}{2} + \gamma_t \frac{2 \bar{\rho}(C_g + L C_p)}{t} \eqs.
\eeq
As $\gamma_t = 2/(t+1)$, from a high-level point of view, the above inequality behaves similarly to
$h_{t+1} \leq (1-(1/t))h_t + {\cal O}(1/t^2)$.
Consequently, applying \cite[Lemma~4]{polyak87} yields a ${\cal O}(1/t)$ convergence rate for $h_t$.
In fact, this is a deterministic version of the case analyzed by \cite[Theorem~10]{ourpaper}. In particular, setting $\alpha = 1, K = 2$
in \cite[(56)]{ourpaper} and using an induction argument yield
\beq
h_t \leq 2 \cdot (4\bar{\rho}(C_g + L C_p) + L \bar{\rho}^2) / (t+1),~\forall~t \geq 1 \eqs.
\eeq

\noindent \textbf{Case 2}: For the case when $F$ is $\mu$-strongly convex and
$\prm^\star$ lies in the interior of $\Cset$ with distance $\delta > 0$ (cf.~\eqref{eq:int}).
Using \cite[Lemma 6]{ourpaper}, we have
\beq
\langle \bar{\prm}_t - \bar{\atom}_t, \grd F(\bar{\prm}_t) \rangle \geq \sqrt{2 \mu \delta^2 h_t } \eqs.
\eeq
Plugging the above into \eqref{eq:fw_key2} gives
\beq
h_{t+1} \leq \sqrt{h_t} ( \sqrt{h_t} - \gamma_t \sqrt{2 \mu \delta^2} ) + \gamma_t^2 \frac{L \bar{\rho}^2}{2} + \gamma_t \frac{2\bar{\rho}(C_g + L C_p)}{t} \eqs.
\eeq
Compared to the case analyzed in \eqref{eq:slow_conv}, when $h_t$ is decreased,
the decrement in $h_{t+1}$ is increased, leading to a faster convergence.
This is a deterministic version of the case analyzed in \cite[Theorem~7]{ourpaper}. Setting $\alpha = 1, K=2$ in \cite[(48)]{ourpaper}
and using an induction argument yields
\beq
h_t \leq  \frac{ (4\bar{\rho}(C_g + L C_p ) + L \bar{\rho}^2 )^2 }{2 \delta^2 \mu} \cdot \frac{9}{(t+1)^2},~\forall~t \geq 1 \eqs.
\eeq

\section{Proof of Theorem~\ref{thm:ncvx}} \label{pf:ncvx}

\subsection{Convergence rate}
Let us recall the definition of the \emph{FW gap}:
\beq
g_t \eqdef \max_{ \prm \in \Cset }~\langle \grd F ( \bar{\prm}_t ), \bar{\prm}_t - \prm \rangle = \langle \grd F ( \bar{\prm}_t ), \bar{\prm}_t - \bar{\atom}_t \rangle \eqs,
\eeq
where we have used the definition of $\bar{\atom}_t$ 
in \eqref{eq:abar} from the previous proof.
For simplicity, we shall assume that $T$ is an even integer in the following.

From the $L$-smoothness of $F$, we have:
\beq \label{eq:key_ineq}
F( \bar{\prm}_{t+1} ) \leq F( \bar{\prm}_t ) + \langle \grd F( \bar{\prm}_t ) , \bar{\prm}_{t+1} - \bar{\prm}_t \rangle + \frac{L}{2} \| \bar{\prm}_{t+1} - \bar{\prm}_t \|_2^2 \eqs.
\eeq
Observe that:
\beq \textstyle
\bar{\prm}_{t+1} - \bar{\prm}_t  = N^{-1} \sum_{i=1}^N \gamma_t (\atom_i^{t} - \bar{\prm}_t^i ) \eqs.
\eeq
As $\atom_t^i$, $\bar{\prm}_t^i \in \Cset$, we have $\| \bar{\prm}_{t+1} - \bar{\prm}_t \|_2 \leq \gamma_t \bar{\rho}$.
Using \eqref{eq:err_bd} and \eqref{eq:thm1_result} from the previous proof of Theorem~\ref{thm:cvx},
the inequality \eqref{eq:key_ineq} can be bounded as:
\beq \label{eq:ineq}
\begin{split}
F( \bar{\prm}_{t+1} ) \leq &~ \textstyle F( \bar{\prm}_t ) - \gamma_t \langle \grd F(\bar{\prm}_t), \bar{\prm}_t - \bar{\atom}_t \rangle \\
& + 2 \gamma_t \bar\rho \cdot ( \dgrd_t + L \cdot \dprm_t ) + \gamma_t^2 L \bar{\rho}^2 / 2 \\
& \hspace{-1.7cm} = F( \bar{\prm}_t) - \gamma_t  g_t +  2 \gamma_t \bar\rho \cdot ( \dgrd_t + L \cdot \dprm_t ) + \gamma_t^2 \frac{ L \bar{\rho}^2 }{ 2 } \eqs.
\end{split}
\eeq
From the definition, we observe that $g_t \geq 0$. Now, summing the two sides of \eqref{eq:ineq} from
$t=T/2+1$ to $t=T$ gives:
\beq
\begin{split}
 \sum_{t=T/2+1}^T \gamma_t g_t ~&  \leq \sum_{t=T/2+1}^T \Big( F( \bar{\prm}_t ) - F (\bar{\prm}_{t+1} ) \Big) \\
&  \hspace{-1.5cm} + \sum_{t=T/2+1}^T \Big(2 \gamma_t \bar{\rho} \cdot ( \dgrd_t + L \cdot \dprm_t ) + \gamma_t^2 \frac{L \bar{\rho}^2 }{2} \Big) \eqs.
\end{split}
\eeq
Canceling duplicated terms in the first term of the right hand side above gives:
\beq \label{eq:pure}
\begin{split}
 \sum_{t=T/2+1}^T \gamma_t g_t ~\leq & ~F(\bar{\prm}_{T/2+1}) - F(\bar{\prm}_{T+1}) \\
&  \hspace{-2cm} + \sum_{t=T/2+1}^T \big(2 \gamma_t \bar{\rho} \cdot ( \dgrd_t + L \cdot \dprm_t ) + \gamma_t^2 \frac{\Curve}{2} \big) \eqs.
\end{split}
\eeq
As $g_t, \gamma_t \geq 0$, we can lower bound the left hand side as:
\beq
 \sum_{t=T/2+1}^T \gamma_t g_t  \geq \Big(\min_{ t \in [T/2+1,T] } g_t \Big) \cdot \Big(\sum_{t=T/2+1}^T \gamma_t \Big) \eqs,
\eeq
and observe that for all $T \geq 6$ and $\alpha \in (0,1)$,
\beq \label{eq:gam_lb}
 \sum_{t=T/2+1}^T \gamma_t \geq \frac{T ^{1-\alpha}}{1-\alpha} \cdot \Big( 1 - \Big(\frac{2}{3} \Big)^{1-\alpha} \Big) = \Omega(T^{1-\alpha}) \eqs.
\eeq
Define the constant
$C \eqdef L \bar{\rho}^2 / 2 + 2 \bar{\rho} ( C_g + L C_p)$.
When $\alpha \geq 0.5$, using the fact that $\gamma_t = t^{-\alpha}$, $\dprm_t = C_p / t^\alpha$, $\dgrd_t = C_g / t^\alpha$,
the right hand side of \eqref{eq:pure} is bounded above by:
\beq \textstyle
G \cdot \rho + C \cdot \sum_{t=T/2+1}^T t^{-2\alpha} \leq G \cdot \rho + C \cdot \log 2 \eqs,
\eeq
note that the series is converging as we are summing from $t=T/2+1$ to $t=T$.
Dividing the above term by the lower bound \eqref{eq:gam_lb}
to $\sum_{t=T/2+1}^T \gamma_t$ yields \eqref{eq:mainthm}.

On the other hand, when $\alpha < 0.5$, we notice that
\beq
\sum_{t=T/2+1}^T t^{-2\alpha} \leq \int_{T/2}^T t^{-2\alpha}~dt = \frac{ 2^{1-2\alpha} - 1}{1 - 2 \alpha} \Big( \frac{T}{2} \Big)^{1-2\alpha} .
\eeq
Therefore, the right hand side of
\eqref{eq:pure} is bounded above by
\beq
G \rho + C \sum_{t=T/2+1}^T t^{-2\alpha}
\leq \Big( G \rho + C \frac{1 - (1/2)^{1-2\alpha}}{1-2\alpha} \Big) \cdot T^{1-2\alpha}  \eqs.
\eeq
Dividing the above term by the lower bound \eqref{eq:gam_lb}
to $\sum_{t=T/2+1}^T \gamma_t$ yields \eqref{eq:mainthm2}.

\subsection{Convergence to stationary point of \eqref{eq:opt}}
Recall that the set of stationary points to \eqref{eq:opt} is defined as:
\beq \textstyle
\Cset^\star \eqdef \{ \bar\prm \in \Cset : \max_{\prm \in \Cset} \langle \grd F( \bar\prm ), \bar\prm - \prm \rangle = 0 \}  \eqs.
\eeq
We state the following Nurminskii's sufficient condition:
\begin{Theorem} \label{thm:nurmin} \cite[Theorem 1]{nurminskii72}
Consider a sequence $\{ \bar{\prm}_t \}_{t \geq 1}$ in a compact set $\Cset$.
Suppose that the following hold\footnote{To give a
clearer presentation, we have rephrased
conditions A.2 and A.3 from the original Nurminskii's conditions.}:
\begin{enumerate}
\item[A.1] $\lim_{t \rightarrow \infty} \| \bar{\prm}_{t+1} - \bar{\prm}_{t} \| = 0$.
\item[A.2] Let $\underline{\prm}$ be a limit point of $\{ \bar{\prm}_t \}_{t \geq 1}$ and
$\{ \prm_{s_t} \}_{t \geq 1}$ be a subsequence that converges to $\underline{\prm}$.
If $\underline{\prm} \notin \Cset^\star$, then
for any $t$ and some sufficiently small  $\epsilon > 0$,
there exists a finite $s$ such that $\| \bar{\prm}_s - \bar{\prm}_{s_t} \| > \epsilon$ and $s > s_t$.
\item[A.3] Let $\underline{\prm}$ be a limit point of $\{ \bar{\prm}_t \}_{t \geq 1}$ and
$\{ \prm_{s_t} \}_{t \geq 1}$ be a subsequence that converges to $\underline{\prm}$.
If $\underline{\prm} \notin \Cset^\star$, then for any $t$ and some sufficiently small
$\epsilon > 0$, we can define
\beq \label{eq:condition2}
\tau_t \eqdef \min_{ s > s_t }~ s~~{\rm s.t.}~~\| \bar{\prm}_s - \bar{\prm}_{s_t} \| > \epsilon
\eeq
where $\tau_t$ is finite. Also,
there exists a continuous function $W(\bar{\prm})$
that takes a finite number of values in $\Cset^\star$
with
\beq \label{eq:lastcond}
\limsup\limits_{t \rightarrow \infty}~ W( \bar{\prm}_{\tau_t} ) < \lim_{t \rightarrow \infty} W( \bar{\prm}_{s_t} ) \eqs.
\eeq
\end{enumerate}
Then the sequence $\{ W( \bar{\prm}_t ) \}_{t \geq 1}$ converges and the limit points of the sequence $\{ \bar{\prm}_t \}_{t \geq 1}$ belongs to the set $\Cset^\star$.
\end{Theorem}

We apply the above theorem to prove that every limit point of
$\{ \bar{\prm}_t \}_{t \geq 1}$ are in $\Cset^\star$.
First, A.1 can be easily verified since
\beq
\| \bar{\prm}_{t+1} - \bar{\prm}_{t} \| \leq \frac{\gamma_{t}}{N} \sum_{i=1}^N \| \atom_{t}^i - \bar{\prm}_{t} \| \leq \frac{ \gamma_t \bar{\rho} }{N} \eqs
\eeq
and we have $\gamma_{t} \rightarrow 0$ as $t \rightarrow \infty$.

As $\Cset$ is compact, there
exists a convergent subsequence $\{ \bar{\prm}_{s_t} \}_{t \geq 1}$ of the
sequence of iterates generated by
the DeFW algorithm.
Let $\underline{\prm}$ be the limit point of $\{ \bar{\prm}_{s_t} \}_{t \geq 1}$
and $\underline{\prm} \notin \Cset^\star$. We shall verify A.2 by contradiction.
In particular, fix a sufficiently small $\epsilon > 0$
and assume that the following holds:
\beq \label{eq:contradict}
\| \bar{\prm}_s - \bar{\prm}_{s_t} \| \leq \epsilon,~\forall~s > s_t,~\forall~t \geq 1 \eqs.
\eeq
As $\{ \bar{\prm}_{s_t} \}_{t \geq 1}$ converges to
$\underline{\prm}$, the assumption \eqref{eq:contradict}
implies that for some sufficiently large $t$ and any $s > s_t$,
we have $\bar{\prm}_s \in {\cal B}_{2 \epsilon} ( \underline{\prm} )$, i.e.,
the ball of radius $2 \epsilon$ centered at $\underline{\prm}$.

Since $\underline{\prm} \notin \Cset^\star$,
the following holds for some $\delta >0$,
\beq \label{eq:contrad}
\langle \grd F( \bar{\prm}_s ), \prm - \bar{\prm}_s \rangle \leq -\delta < 0,~\forall~\prm \in \Cset,~\forall s > s_t \eqs.
\eeq
In particular, we have
$\langle \grd F( \bar{\prm}_s ), \bar\atom_s - \bar{\prm}_s \rangle \leq -\delta$
as we recall that
$\bar\atom_s = \arg \min_{ \atom \in \Cset } \langle \grd F( \bar{\prm}_s ), \atom \rangle$.

On the other hand, from \eqref{eq:ineq} and using Assumptions~\ref{ass:a1}-\ref{ass:a2},
it holds true for all $t \geq 1$ that:
\beq \label{eq:someineq}
\begin{split}
F( \bar{\prm}_{t+1} ) - F( \bar{\prm}_t ) & \leq \gamma_t \cdot \langle \grd F( \bar{\prm}_t ), \bar\atom_t - \bar{\prm}_t \rangle\\
& \hspace{.3cm} + \gamma_t \cdot {\cal O} ( t^{-\alpha} ) + \gamma_t^2 L \bar{\rho}^2 / 2 \eqs.
\end{split}
\eeq
To arrive at a contradiction, we let $s > s_t$ and sum up the two sides of \eqref{eq:someineq}
from $t=s_t$ to $t=s$. Consider the following chain of inequality:
\beq \label{eq:tauk}
\begin{split}
F( \bar{\prm}_s ) - F( \bar{\prm}_{s_t} ) & \leq \sum_{\ell = s_t}^s \gamma_\ell ( \grd F( \bar{\prm}_\ell), \bar\atom_\ell - \bar{\prm}_\ell \rangle + {\cal O}(\ell^{-\alpha}) ) \\
&  \leq - \delta \sum_{\ell=s_t}^s \gamma_\ell + \sum_{\ell=s_t}^s \gamma_\ell {\cal O}( \ell^{-\alpha} ) \eqs,
\ifplainver 
\\
\else
\\[-.6cm]
\fi
\end{split} \vspace{.4cm}
\eeq
where the first inequality is due to the fact that
$\gamma_\ell^2 L \bar{\rho}^2 / 2 = \gamma_\ell {\cal O}(\ell^{-\alpha})$
and the second inequality is due to \eqref{eq:contrad}.
Rearranging terms in \eqref{eq:tauk}, we have
\beq \label{eq:infinity}
F( \bar{\prm}_s ) - F( \bar{\prm}_{s_t} ) - \sum_{\ell=s_t}^s C \cdot \ell^{-2\alpha}
\leq  - \delta \sum_{\ell=s_t}^s \ell^{-\alpha} \eqs,
\eeq
for some $C < \infty$.
As $1 \geq \alpha > 0.5$, we have $\lim_{s \rightarrow \infty} \sum_{\ell=s_t}^s \ell^{-2\alpha} < \infty$
on the left hand side and $\lim_{s \rightarrow \infty} \sum_{\ell=s_t}^s \ell^{-\alpha} \rightarrow + \infty$ on the right hand side.
Letting $s \rightarrow \infty$ on the both side of \eqref{eq:infinity}
implies
\beq
\lim_{s \rightarrow \infty} F( \bar{\prm}_s ) - F( \bar{\prm}_{s_t} ) < - \infty \eqs,
\eeq
This leads to a contradiction to \eqref{eq:contrad} since $F (\prm)$ is bounded over $\Cset$.
We conclude that A.2 holds for the DeFW algorithm.

The remaining task is to verify A.3. We notice that the indices $\tau_t$ in \eqref{eq:condition2}
are well defined since A.2 holds. Take $W( \prm ) = F( \prm )$ and notice
that the image $F( \Cset^\star )$ is a finite set (cf.~Assumption~\ref{ass:a3}).
By the definition of $\tau_t$,
we have $\bar{\prm}_{s} \in {\cal B}_\epsilon ( \bar{\prm}_{s_t} )$ for all $s_t \leq s \leq \tau_t - 1$.
Again for some sufficiently large $t$, we have $\bar{\prm}_{s} \in {\cal B}_\epsilon ( \bar{\prm}_{s_t} ) \subseteq {\cal B}_{2 \epsilon} ( \underline{\prm} )$ and the inequality \eqref{eq:tauk} holds
for $s= \tau_t - 1$. This gives:
\beq \label{eq:thirdcond}
F( \bar{\prm}_{\tau_t} ) - F (\bar{\prm}_{s_t} ) \leq \sum_{\ell= s_t}^{\tau_t - 1} \gamma_\ell \cdot ( -\delta + {\cal O}( \ell^{-\alpha} ) ) \eqs.
\eeq
On the other hand, we have $\bar{\prm}_{\tau_t} \notin {\cal B}_\epsilon ( \bar{\prm}_{s_t} )$ and thus
\beq
\epsilon < \| \bar{\prm}_{\tau_t} - \bar{\prm}_{s_t} \| \leq \sum_{\ell=s_t}^{\tau_t - 1} \gamma_\ell \Big\| \sum_{i=1}^N \frac{\atom_\ell^i}{N} - \bar{\prm}_\ell \Big\| \leq \bar{\rho} \sum_{\ell=s_t}^{\tau_t - 1} \gamma_\ell \eqs.
\eeq
The above implies that $\sum_{\ell= s_t}^{\tau_t - 1} \gamma_\ell > \epsilon / \bar{\rho} > 0$.
Considering \eqref{eq:thirdcond} again, observe that ${\cal O}( \ell^{-\alpha} )$ decays to zero,
for some sufficiently large $t$, we have $- \delta + {\cal O} (\ell^{-\alpha} ) \leq -\delta' < 0$
if $\ell \geq s_t$. Therefore, \eqref{eq:thirdcond} leads to
\beq
F( \bar{\prm}_{\tau_t} ) - F (\bar{\prm}_{s_t} ) \leq - \delta' \sum_{\ell=s_t}^{\tau_t - 1} \gamma_\ell < - \frac{ \delta' \epsilon } {\bar{\rho} } < 0 \eqs.
\eeq
Taking the limit $t \rightarrow \infty$ on both sides leads to \eqref{eq:lastcond}.
The proof for the convergence to stationary point
in Theorem~\ref{thm:ncvx} is completed by applying Theorem~\ref{thm:nurmin}.

\section{Proof of Lemma~\ref{lem:gac2}} \label{pf:gac2}
For simplicity, we shall drop the dependence of $\alpha$ in the constant $t_0(\alpha)$.
It suffices to show that
for all $t \geq 1$,
\beq \label{eq:ind}
\sqrt{ \sum_{i=1}^N \| \bar{\prm}_t^i - \bar{\prm}_t \|_2^2 } \leq \frac{C_p}{ t^\alpha },~C_p = (t_0)^\alpha \cdot \sqrt{N} \bar{\rho} \eqs.
\eeq
We observe that for $t = 1$ to $t = t_0$, the above inequality is true
since $\bar{\prm}_t^i, \bar{\prm}_t \in \Cset$ and the
diameter of $\Cset$ is bounded by $\bar{\rho}$.
For the induction step, let us assume that
$\sqrt{ \sum_{i=1}^N \| \bar{\prm}_t^i - \bar{\prm}_t \|^2 } \leq C_p / t^\alpha$ for some $t \geq t_0$. Observe that
\beq
\prm_{t+1}^i = (1- t^{-\alpha}) \bar{\prm}_t^i + t^{-\alpha} \atom_t^i \eqs.
\eeq
Denote $\tilde{\atom}_t = N^{-1} \sum_{i=1}^N \atom_t^i$ and 
using Fact~\ref{fact:gac}, we observe that,
\beq  \begin{split}
{\sum_{i=1}^N \| \bar{\prm}_{t+1}^i - \bar{\prm}_{t+1} \|_2^2 } \leq & \\
& \hspace{-3cm} |\lambda_2 ({\bm W})|^2 \cdot {\sum_{j=1}^N \| (1-t^{-\alpha}) (\bar{\prm}_t^j - \bar{\prm}_t)
+ t^{-\alpha} (\atom_t^j - \tilde{\atom}_t) \|_2^2} \eqs,
\end{split}
\eeq
where we have used the fact $\bar{\prm}_{t+1} =
(1-t^{-\alpha}) \bar{\prm}_t + t^{-\alpha} \tilde{\atom}_t$. The right hand side
in the above can be bounded by
\beq \begin{split}
{\sum_{j=1}^N \| (1-t^{-\alpha}) (\bar{\prm}_t^j - \bar{\prm}_t)
+ t^{-\alpha} (\atom_t^j - \tilde{\atom}_t) \|_2^2} & \\
& \hspace{-6cm} \leq \sum_{j=1}^N \big(  \| \bar{\prm}_t^j - \bar{\prm}_t \|_2^2 + t^{-2\alpha} \bar{\rho}^2 + 2 \bar{\rho} t^{-\alpha} \| \bar{\prm}_t^j - \bar{\prm}_t \|_2 \big) \\
& \hspace{-6cm} \leq t^{-2 \alpha} ( C_p^2 + N \bar{\rho}^2 ) + 2 \bar{\rho} t^{-\alpha} \sqrt{N} \sqrt{\sum_{j=1}^N \| \bar{\prm}_t^j - \bar{\prm}_t \|_2^2 }  \\
& \hspace{-6cm} \leq t^{-2 \alpha} ( C_p + \sqrt{N} \bar{\rho})^2 \leq \Big( \frac{(t_0)^\alpha+1}{(t_0)^\alpha \cdot t^{\alpha}} \cdot C_p \Big)^2 \eqs,
\end{split}
\eeq
where we have used the boundedness of $\Cset$ in the first inequality,
the norm equivalence $\sum_{j=1}^N |c_j| \leq \sqrt{N} \sqrt{\sum_{j=1}^N c_j^2}$ in the second inequality
and the induction hypothesis in the third and fourth inequalities.
Consequently, from \eqref{eq:t0}, we observe that for all $t \geq t_0$,
\beq
|\lambda_2 ({\bm W})| \cdot \frac{(t_0)^\alpha+1}{(t_0)^\alpha \cdot t^{\alpha}} \leq \frac{1}{(t+1)^{\alpha}} \eqs,
\eeq
and the induction step is completed.
Finally, Lemma~\ref{lem:gac2} is proven by observing that \eqref{eq:ind} implies \eqref{eq:gac2r}.

\section{Proof of Lemma~\ref{lem:gac1}} \label{pf:gac1}
We prove the first condition \eqref{eq:avg1} with a simple induction.
This condition is obviously true for the base step $t=1$. For induction step, suppose that
\eqref{eq:avg1} is true up to some $t$, then
\beq \begin{split}
 \sum_{i=1}^N \bgrd{t+1}{i} = & \sum_{i=1}^N (\bargrd{t}{i} - \grd f_i (\bar{\prm}_t^i)) + \sum_{i=1}^N \grd f_i ( \bar{\prm}_{t+1}^i ) \eqs.
\end{split}
\eeq
Note that the first term on the right hand side is zero due to the induction hypothesis.
Thus, the induction step is completed and
$N^{-1} \sum_{i=1}^N \bgrd{t}{i} = N^{-1} \sum_{i=1}^N \grd f_i (\bar{\prm}_t^i )$ for all $t \geq 1$.
Lastly, as ${\bm W}$ is doubly stochastic, we have 
$N^{-1} \sum_{i=1}^N \bargrd{t}{i} = N^{-1} \sum_{i=1}^N \bgrd{t}{i}$.

Then, we prove the second condition \eqref{eq:gac1r}.
For simplicity, we drop the dependence of $\alpha$ in the constant $t_0(\alpha)$.
Recall $\bargrd{t}{} \eqdef N^{-1} \sum_{i=1}^N \bgrd{t}{i}$.
It suffices to prove:
\beq \label{eq:ind2}
\begin{split}
& \sqrt{ \sum_{i=1}^N \| \bargrd{t}{i} - \bargrd{t}{} \|_2^2 } \leq \frac{C_g}{t^\alpha} \eqs, \\
& C_g = \sqrt{N} \max \Big\{ |\lambda_2({\bm W})| (t_0)^\alpha \Big( \frac{ L \bar{\rho} }{1 - |\lambda_2({\bm W})|}  + B_1 \Big), 2 (2C_p + \bar{\rho} ) L \Big\} \eqs,
\end{split}
\eeq
for all $t \geq 1$ using induction.
For $t=1$ to $t = t_0$, we shall prove that the left hand side of the inequality
is bounded. To proceed,
we define the $d \times N$ matrices:
\beq
\begin{split}
& {\bm E}_t \eqdef \big( \bargrd{t}{1} \cdots \bargrd{t}{N} \big) - \big(\bargrd{t}{} \cdots 
\bargrd{t}{} \big) \\
& \grd {\bm F}_t \eqdef ( \grd f_1( \bar{\prm}_t^1 ) \cdots \grd f_N( \bar{\prm}_t^N ) ) \eqs,
\end{split}
\eeq
and observe
\beq \textstyle
\sum_{i=1}^N \| \bargrd{t}{i} - \bargrd{t}{} \|_2^2= \| {\rm vec}( {\bm E}_t ) \|_2^2 \eqs.
\eeq
Furthermore, ${\bm E}_1 = \grd {\bm F}_1 ( {\bm W} - (1/N) {\bf 1}{\bf 1}^\top )$, 
and we have the following recursion for $t \geq 2$, 
\beq
\begin{split}
{\bm E}_t & = \big( \bargrd{t}{1} \cdots \bargrd{t}{N} \big) - \big(\bargrd{t}{} \cdots 
\bargrd{t}{} \big) \\
& = \Big( \big( \bargrd{t-1}{1} \cdots \bargrd{t-1}{N} \big) + \grd {\bm F}_t - \grd {\bm F}_{t-1} \Big) {\bm W} - \big(\bargrd{t}{} \cdots 
\bargrd{t}{} \big) \\
& = ( {\bm E}_{t-1} + \grd {\bm F}_t - \grd {\bm F}_{t-1} ) ( {\bm W} - (1/N) {\bf 1}{\bf 1}^\top ) \eqs,
\end{split}
\eeq 
where we have used the equivalence below:
\beq
\frac{1}{N} \grd {\bm F}_t {\bf 1}{\bf1}^\top = \frac{1}{N} \big( \bargrd{t}{1} \cdots \bargrd{t}{N} \big) {\bf 1}{\bf 1}^\top =  \big(\bargrd{t}{} \cdots \bargrd{t}{} \big) \eqs.
\eeq

For $t = 1$,  
$\| {\rm vec} ( {\bm E}_1 ) \|_2 \leq |\lambda_2({\bm W})| \sqrt{N} B$
since $\| {\rm vec}( \grd {\bm F}_1 ) \|_2 \leq \sqrt{N} B_1$.
For $t \geq 2$, we have
\beq
\begin{split}
& \| {\rm vec} ( {\bm E}_t ) \|_2 \leq \| ( {\bm W} - (1/N) {\bf 1}{\bf 1}^\top ) \otimes {\bm I} \|_2  \\
& \hspace{.2cm} \big( \| {\rm vec} ( {\bm E}_{t-1} ) \|_2 + \| {\rm vec} ( \grd {\bm F}_t  - \grd {\bm F}_{t-1} ) \|_2 \big) \eqs.
\end{split}
\eeq
Since
$\| ( {\bm W} - (1/N) {\bf 1}{\bf 1}^\top ) 
\otimes {\bm I} \|_2 \leq |\lambda_2( {\bm W})|$
and the $L$-smoothness of $f_i$ implies
$\| {\rm vec} ( \grd {\bm F}_t  - \grd {\bm F}_{t-1} ) \|_2 \leq \sqrt{N} L \bar{\rho}$,
this leads to
\beq
\| {\rm vec} ( {\bm E}_t ) \|_2 \leq 
|\lambda_2( {\bm W})| \cdot \big( \| {\rm vec} ( {\bm E}_{t-1} ) \|_2 + \sqrt{N} L \bar{\rho} \big) 
\eeq
Evaluating the recursion above shows 
\beq
\| {\rm vec}({\bm E}_t) \|_2 \leq \frac{ |\lambda_2({\bm W})| \sqrt{N} L \bar{\rho} }{1 - |\lambda_2({\bm W})|}  \!~+ |\lambda_2({\bm W})| \sqrt{N} B_1 \eqs,
\eeq 
and thus proving the base step for $t=1$ to $t=t_0$. 

For the induction step with $t > t_0$, we suppose that
$\sqrt{ \sum_{i=1}^N  \| \bargrd{t}{i} - \bargrd{t}{} \|_2^2 } \leq C_g / t^\alpha$ for some $t \geq t_0$.
Define the slack variable $\delta f_{t+1}^i \eqdef \grd f_i(\bar{\prm}_{t+1}^i) - \grd f_i( \bar{\prm}_t^i )$.
{We observe that $\bgrd{t+1}{i} = \delta f_{t+1}^i + \bargrd{t}{i}$
and $\bargrd{t+1}{i} = \sum_{j=1}^N W_{ij} \bgrd{t+1}{j}$, thus
applying Fact~\ref{fact:gac} yields}
\beq \label{eq:lem2_1st}
\begin{split}
\sum_{i=1}^N  \| \bargrd{t+1}{i} - \bargrd{t+1}{} \|_2^2 & \leq \\
& \hspace{-3cm} |\lambda_2 ({\bm W})|^2 \cdot \sum_{i=1}^N \| \bargrd{t}{i} + \delta f_{t+1}^i -  \bargrd{t+1}{} \|_2^2 \eqs,
\end{split}
\eeq
Similarly, define $\delta F_{t+1} \eqdef \bargrd{t+1}{} - \bargrd{t}{} =  N^{-1}\sum_{i=1}^N \delta f_{t+1}^i$
and observe that
we can bound the right hand side of \eqref{eq:lem2_1st} as
\beq \label{eq:longlem} \begin{split}
\sum_{i=1}^N \| \bargrd{t}{i} + \delta f_{t+1}^i - \bargrd{t+1}{} \|_2^2 & \\
& \hspace{-4.5cm} \leq \sum_{i=1}^N \Big( \| \bargrd{t}{i} - \bargrd{t}{} \|_2^2 + \| \delta f_{t+1}^i - \delta F_{t+1} \|_2^2  \\
& \hspace{-3.6cm} + 2 \cdot \| \delta f_{t+1}^i - \delta F_{t+1} \|_2 \cdot \| \bargrd{t}{i} - \bargrd{t}{} \|_2 \Big)
\end{split}
\eeq
where the first inequality is obtained by expanding the squared $\ell_2$ norm and applying
Cauchy-Schwartz inequality.

Observe that for all $i \in [N]$, we have the following chain:
\beq 
\begin{split}
\| \delta f_{t+1}^i \|_2 = \| \grd f_i (\bar{\prm}_{t+1}^i) - \grd f_i ( \bar{\prm}_t^i ) \|_2 & \leq L \| \bar{\prm}_{t+1}^i - \bar{\prm}_t^i \|_2 \\
& \hspace{-5.3cm} \textstyle \leq L \Big\| \sum_{j=1}^N W_{ij} \big( (\prm_{t+1}^j - \bar{\prm}_t^j ) + (\bar{\prm}_t^j - \bar{\prm}_t^i) \big) \Big\|_2 \\
& \hspace{-5.3cm} \textstyle \leq L \sum_{j=1}^N W_{ij} \Big( t^{-\alpha} \bar{\rho} + 2 C_p t^{-\alpha} \Big) = ( 2 C_p + \bar{\rho} ) L t^{-\alpha} \eqs, \vspace{2cm}
\end{split}\vspace{-2cm}
\eeq
where the last inequality is due to the convexity of $\ell_2$ norm, the update rule in line~\ref{line:fw} of Algorithm~\ref{alg:defw} and the results from Lemma~\ref{lem:gac2}.
Using the triangular inequality, we observe that
\beq
\begin{split}
\| \delta f_{t+1}^i  - \delta F_{t+1} \|_2 & \textstyle = \Big\| \big(1 - \frac{1}{N} \big) \delta_{t+1}^i + \frac{1}{N} \sum_{j \neq i} \delta_{t+1}^j \Big\|_2 \\
& \hspace{-2.3cm} \textstyle \leq \big(1 - \frac{1}{N} \big) \| \delta_{t+1}^i \|_2 + \frac{1}{N} \sum_{j \neq i} \| \delta_{t+1}^j \|_2 \\
& \hspace{-2.3cm} \leq 2 \big(1 - \frac{1}{N} \big) ( 2 C_p + \bar{\rho} ) L t^{-\alpha} \leq 2 ( 2 C_p + \bar{\rho} ) L t^{-\alpha} \eqs.
\end{split}
\eeq
Finally, applying the induction hypothesis, the right hand side of
Eq.~\eqref{eq:longlem} can be bounded by
\beq \notag
 \begin{split}
\sum_{i=1}^N \| \bargrd{t}{i} + \delta f_{t+1}^i - \bargrd{t+1}{} \|_2^2 & \\
& \hspace{-4.4cm} \leq t^{-2\alpha} \big( C_g^2 + 4N (2C_p + \bar{\rho} )^2 L^2 \big) \\
& \hspace{-3.8cm} \textstyle + t^{-\alpha} 4L ( 2 C_p + \bar{\rho} )  \sqrt{N} \sqrt{ \sum_{i=1}^N   \| \bargrd{t}{i} - \bargrd{t}{} \|_2^2} \\[.1cm]
& \hspace{-4.4cm} \leq t^{-2\alpha} \!~ \big( C_g + 2L \sqrt{N}(2 C_p + \bar{\rho} ) \big)^2 \leq \Big( \frac{(t_0)^\alpha+1}{(t_0)^\alpha \cdot t^{\alpha}} \cdot C_g \Big)^2,
\end{split}
\eeq
where we have used the fact that $\sum_{i=1}^N \| \bargrd{t}{i} - \bargrd{t}{} \|_2 \leq \sqrt{N}
\sqrt{ \sum_{i=1}^N \| \bargrd{t}{i} - \bargrd{t}{} \|_2^2 }$ in the first inequality.
Invoking \eqref{eq:t0}, we can upper bound the right hand side of \eqref{eq:lem2_1st} by $C_g^2 / (t+1)^{2\alpha}$
for all $t \geq t_0$.
Taking square root on both sides of the inequality completes the induction step.
Consequently, \eqref{eq:gac1r} can be implied by \eqref{eq:ind2}.

\section{Proof of Lemma~\ref{lem:prob}} \label{pf:prob}
Applying triangular inequality on the error vector yields:
\begin{equation} \label{eq:probab}
\begin{split}
 \Big\| \xi_t^{-1} \bargrd{t}{i} - \frac{1}{N} \sum_{j=1}^N \grd f_j (\bar{\prm}_t^j ) \Big\|_\infty & \\
& \hspace{-5cm} \leq  \xi_t^{-1} \cdot \Big\| \bargrd{t}{i} -  \frac{1}{N} \sum_{j=1}^N \grd f_j (\bar{\prm}_t^j ) \odot {\bf 1}_{\Omega_t} \Big\|_\infty \\
& \hspace{-5cm}  ~~ + \Big\| \big( \frac{1}{N} \sum_{j=1}^N \grd f_j (\bar{\prm}_t^j ) \big) \odot (\xi_t^{-1} {\bf 1}_{\Omega_t}  - {\bm 1})  \Big\|_\infty \eqs,
\end{split}
\end{equation}
where ${\bf 1}$ denotes the all-one vector.
For the first term in the right hand side of \eqref{eq:probab}, observe that $\bargrd{t}{i}$ is obtained by applying the AC updates
on the sparsified local gradients $\grd f_i( \bar{\prm}_t^i) \odot {\bf 1}_{\Omega_t}$ for
$\ell_t =  \lceil C_l + \log ( t ) / \log |\lambda_2^{-1} ( {\bm W} )| \rceil$ rounds, applying Fact~\ref{fact:gac} yields the following for
all $i \in [N]$:
\beq
\begin{split}
 \Big\| \bargrd{t}{i} -  \frac{1}{N} \sum_{j=1}^N \grd f_j (\bar{\prm}_t^j ) \odot {\bf 1}_{\Omega_t} \Big\|_\infty & \\
& \hspace{-5.2cm} \leq  | \lambda_2 ( {\bm W} )|^{\ell_t} \cdot \Big\| ( \grd f_i (\bar{\prm}_t^i ) - \frac{1}{N} \sum_{j=1}^N \grd f_j (\bar{\prm}_t^j ) ) \odot {\bf 1}_{\Omega_t} \Big\|_\infty \\
& \hspace{-5.2cm} \leq  |\lambda_2 ( {\bm W} )|^{C_l} \cdot B / t \eqs,
\end{split}
\eeq
for some $B < \infty$ since the gradients are bounded.

For the second term in the right hand side of \eqref{eq:probab}, we first apply the inequality $ \| \big( N^{-1} \sum_{i=1}^N \grd f_i (\bar{\prm}_t^i ) \big) \odot (\xi_t^{-1} {\bf 1}_{\Omega_t} - {\bm 1})  \|_\infty \leq  \| N^{-1} \sum_{i=1}^N \grd f_i (\bar{\prm}_t^i ) \|_\infty \| (\xi_t^{-1} {\bf 1}_{\Omega_t} - {\bm 1})  \|_\infty$ from \cite{Horn_Johnson94}.
Now, the probability that coordinate  $k$ is included is given by:
\beq
P ( k \in \Omega_t ) = 1 - P \Big( \bigcap_{i=1}^N k \notin \Omega_{t,i} \Big) = 1 - ( 1 - \frac{1}{d} )^{p_t N} = \xi_t \eqs,
\eeq
and that $\EE[ {\bf 1}_{\Omega_t} ] = \xi_t {\bf 1}$.
Then, observing that each element in $\xi_t^{-1} {\bf 1}_{\Omega_t}$ is bounded
in $[0,\xi_t^{-1}]$ and applying the Hoefding's inequality \cite{pascal03},
the following holds true for all $x > 0$:
\beq
P \big( \| \xi_t^{-1} {\bf 1}_{\Omega_t} - {\bm 1} \|_\infty \geq x \big) \leq 2d \cdot e^{- 2 x^2 / \xi_t^{-2} } \eqs,
\eeq
where we have applied a union bound argument to take care of the $\ell_{\infty}$-norm.

Setting $x = \xi_t^{-1} \sqrt{ (\log(2 dt^2) - \log \epsilon) / 2}$ and applying another union bound
show that with probability at least $1 - (\pi^2 \epsilon / 6)$, the following holds
for all $t \geq 1$:
\begin{equation}
\begin{array}{l}
\displaystyle \Big\| \frac{1}{N} \sum_{i=1}^N \grd f_i (\bar{\prm}_t^i ) \odot (\xi_t^{-1} {\bf 1}_{\Omega_t} - {\bm 1})  \Big\|_\infty \vspace{.1cm} \\
\displaystyle ~\leq \xi_t^{-1} \Big\| \frac{1}{N} \sum_{i=1}^N \grd f_i (\bar{\prm}_t^i ) \Big\|_\infty \sqrt{ \frac{\log(2dt^2 / \epsilon )}{2} } \eqs,
\end{array}
\end{equation}
As $d \gg 0$, we have $\xi_t^{-1} \approx d / (p_t N)$.
Recalling $p_t \geq C_0 t$ yields the desired result in Lemma~\ref{lem:prob}.

\bibliographystyle{IEEEtran}
\bibliography{defw_ref}

\ifplainver 

\else

\begin{IEEEbiography}
[{\includegraphics[width=1in,height=1.25in,clip,keepaspectratio]{./bio/To_2.jpg}}]
{Hoi-To Wai} (S'11) received his B. Eng. (with First Class Honor) and M. Phil. degrees in Electronic Engineering from The Chinese University of Hong Kong (CUHK) in 2010 and 2012, respectively.  He is currently a PhD candidate 
in Electrical Engineering at Arizona State University (ASU), USA. Prior to ASU, he was a PhD student at University of California, Davis, USA from 2013 to 2014. He was a visiting PhD student at T\'el\'ecom ParisTech, Paris, France in summer 2015.

His current research interests are in the broad area of distributed optimization, signal processing, machine learning, with a focus on networked systems including social, biological and sensor networks.
\end{IEEEbiography}

\begin{IEEEbiography}
[{\includegraphics[width=1in,height=1.25in,clip,keepaspectratio]{./bio/jean.jpg}}]
{Jean Lafond} is a third year Ph.D. student
in machine learning at University Paris Saclay. He graduated from ENSAE ParisTech in 2010
and worked in quantitative finance at Nomura London for two years.
He also got a Master Degree in Probability and Statistics from University Paris VII in 2013.
His research mainly focuses on matrix completion methods,
stochastic and online optimization.
\end{IEEEbiography}

\begin{IEEEbiography}
[{\includegraphics[width=1in,height=1.25in,clip,keepaspectratio]{./bio/Anna_2.JPG}}]
{Anna Scaglione} (F'11) is currently a professor in electrical and computer engineering at Arizona State University. She was previously at the University of California at Davis, Cornell University and University of New Mexico.
Her research focuses on various applications of signal processing in network science that include intelligent infrastructure, information systems and social networks.

Dr. Scaglione was elected an IEEE fellow in 2011. She received the 2000 IEEE Signal Processing Transactions Best Paper Award and more recently was honored for the 2013, IEEE Donald G. Fink Prize Paper Award for the best review paper in that year in the IEEE publications, her work with her student earned  2013 IEEE Signal Processing Society Young Author Best Paper Award (Lin Li). She was EIC of the IEEE Signal Processing Letters and served in many other capacities the IEEE Signal Processing and the IEEE Communication societies.
\end{IEEEbiography}

\begin{IEEEbiography}
[{\includegraphics[width=1in,height=1.25in,clip,keepaspectratio]{./bio/eric.JPG}}]
{Eric Moulines} received the Engineering degree from Ecole Polytechnique, Paris, France, in 1984, the Ph.~D.~degree in electrical engineering from Ecole Nationale Sup\'erieure des T\'el\'ecommunications, in 1990. In 1990, he joined the Signal and Image processing department at T\'el\'ecom ParisTech where he became a full professor in 1996. In 2015, he joined the Applied Mathematics Center of Ecole Polytechnique, where he is currently a professor in statistics.

His areas of expertise include computational statistics, machine learning, statistical signal processing and time series analysis. His current research topics cover large-scale (Bayesian) inference with applications to inverse problems and machine learning and non-linear filtering.   He has published more than 100 papers in leading journals of the field. In 1997 and 2006, he receive the Best paper Award of  the IEEE Signal Processing Society  (for  papers in IEEE Trans. On Signal Processing)
He served in the editorial boards of IEEE Trans. On Signal Processing, Signal Processing, Stochastic Processes and Applications, Journal of Statistical Planning and Inference. He was the Editor-in-Chief of Bernoulli from 2013-2016.

E.~Moulines is a EURASIP and Insitute of Mathematical Statistics fellow. He was the recipient of the 2010 Silver Medal from the Centre National de Recherche Scientifique, 2011 Orange prize of the French Academy of Sciences.
\end{IEEEbiography}

\fi

\end{document}